\documentclass[12pt]{article}
\usepackage{graphicx}
\usepackage{amssymb}
\usepackage{amsmath}
\usepackage{amsthm}
\usepackage{cite}
\usepackage{graphicx}
\newcommand{\nn}{\nonumber}
\numberwithin{equation}{section}

\def\BR{{\mathbb R}}
\def\BN{{\mathbb N}}

\def\ov{\overline }

\def\b{{\beta}}
\def\om{{\omega}}
\newtheorem{Pa}{Paper}[section]
\newtheorem{Tm}[Pa]{{\bf Theorem}}
\newtheorem{La}[Pa]{{\bf Lemma}}
\newtheorem{Cy}[Pa]{{\bf Corollary}}
\newtheorem{Rk}[Pa]{{\bf Remark}}

\newtheorem{Ee}[Pa]{{\bf Example}}
\newtheorem{Dn}[Pa]{{\bf Definition}}
\newtheorem{Pn}[Pa]{{\bf Proposition}}

\def\XXint#1#2#3{{\setbox0=\hbox{$#1{#2#3}{\int}$}
     \vcenter{\hbox{$#2#3$}}\kern-.5\wd0}}

\title{Periodic functions: self-intersections, \\ local singular points, and folds}

\author{Lev Sakhnovich}

\date{}

\begin{document}
\maketitle
\begin{abstract}  
The oscillations described by periodic functions  play an important role
in many areas.
In the present paper, we study
periodic functions which belong to the class of the $n$-member chains.
The self-intersection and local singular points of these
periodic
functions are constructed. We consider several classical curves in two and three dimensions.
We also introduce and study two  new classes of periodic
functions: the class of periodic helices and the class of S-torus knots.
 In the last section,
we construct  folds for the $2$-member chains.
For these purposes, we derive and use several results on trigonometric
formulas.
\end{abstract}

\vspace{0.2em}

{\bf Keywords:} Periodic functions, self-intersection points, singular points, rotating
group,  periodic helix, torus knot.

\section{Introduction}
1. Periodic functions describe  oscillations, which
play an important role in many areas (e.g., in mathematics \cite{Bar, Pro},  mechanics \cite{LaLi1}, wave theory  \cite{Fey, LaLi2, SUZ, FKK},
and knot theory \cite{BDF, Col}). Periodic functions $f(t)$ of the form
\begin{equation}f(t)=\sum_{k=1}^{n}c_{k}e^{itm_k}, \quad c_k\in \BR, \quad c_k\not=0 \quad (-\infty<t<\infty),\label{1.1}\end{equation}
where $\BR$ is the real axis and the numbers   $m_k$ are coprime integers, are well-known. The corresponding  curve $f(t)=\phi(t)+i\psi(t)$
$\,\big(\phi(t), \psi(t)\in \BR\big)$ may be written down in the parametric form:
\begin{equation}x=\phi(t)=\sum_{k=1}^{n}c_{k}\cos(tm_k),\quad y=\psi(t)=\sum_{k=1}^{n}c_{k}\sin(tm_k). \label{1.2}\end{equation}
Introduce also a more general case of periodic curves:
\begin{equation}x=\phi(t)=\sum_{k=1}^{n}c_{k}\cos(tm_k),\quad y=\psi(t)=\sum_{k=1}^{n}d_{k}\sin(tm_k) \quad (c_k,d_k\in \BR), \label{1.3}\end{equation}
where $m_k$ are again coprime integers. Formula \eqref{1.3} describes an important subclass of the so called Fourier class.
\emph{We call the expression of the form \eqref{1.3} an $n$-member chain.}

 In the present paper, we study periodic functions which belong to the class \eqref{1.3} of the $n$-member chains.
 \begin{Rk}\label{Rk1} In particular, the paper presents a unified approach for calculating self-intersection points for a wide class of curves.
 The approach is based on our trigonometric formulas $($see section 2$)$ and consists of two steps. First, we find the self-intersection points
lying on the axes $y=0$ or $x=0$ $($see section 3$)$. These points are described in terms of  the roots of
polynomials in one variable. In the second step, we use the results of the first step in order to find general type self-intersection points.
\end{Rk}

2. Let us consider the contents of the paper in greater detail. 
In section 2, we represent $\sin({\ell}x)$ and $\cos({\ell}x)$,  where  $\ell$ belongs to the class $\BN$ of positive integers, as polynomials of $\sin(x)$ and  $\cos(x)$.
In section~3,  we use  trigonometric formulas obtained  in section 2  for  constructing
the self-intersection points of the curves \eqref{1.3}. The important cases where the self-intersection points
belong to the  axes of symmetry  of the curves \eqref{1.3} are studied separately. Section 4 is dedicated to the particular case \eqref{1.1} of the $n$-member chains.
In section 5, we consider $2$-member chains, which belong to the case \eqref{1.1}:
\begin{equation} f(t)=c_{1}e^{itm}+c_{2}e^{it\ell},\quad \label{1.4}\end{equation}
where $m$ and $\ell$ are coprime integers and $\ell>m$.   In this case, we describe the sets of self-intersection points with the help of  groups of rotations (Lie groups).
In section 5, a complete description of all  self-intersection points is given  for the cases
where $\ell-m$ is odd or where $\ell-m=2s$ and $s$ is odd.  We note that  2-member chain  is a main part of the corresponding
periodic function in many  physical problems \cite{SUZ}.  The local singular points  of the curve \eqref{1.3} are constructed in section 6. 
For this purpose, we use the results of section 2.
In section 7, we consider again the 2-member chains \eqref{1.4}. We prove the following assertion.
\begin{Pn}\label{Propsition 1.1}The $2$-member chain \eqref{1.4} has local singular points if and only if
\begin{equation}(mc_1)^{2}=(\ell{c_2})^{2}.\label{1.5}\end{equation}
All the corresponding local singular points are return points of the first kind.\end{Pn}
 A complete description of all local singular points of the 2-member chains \eqref{1.4} is also given in section 7. Similar to the description of the self-intersection
 points in section 5, the sets of singular points are described in section 7  using  groups of rotations of the complex plane (Lie groups).
 
  In section 8, we show that several well known classes of curves (see \cite{Nas, Pro})
 belong to the 2-member chain class: epicycloid, hypocycloid, epitrochoid and hypotrochoid.  In section 9, we use the results of the paper in order to study  3-dimensional curves.
 We introduce a new interesting class of  3-dimensional  curves (periodic helix). Well known 3-dimensional curves (satellite curves, Capareda curves, curve of constant precession \cite{Sco, LSTT} belong to the class of periodic helix curves.
  Section 10 is dedicated to the spectral theory of  integro-differential operators with difference kernels. In this section, we prolong our research on the triangular
 operators with difference kernels
 (see \cite{Sak1,Sak3,Sak4,Sak2}) using results of the present paper. In section 11, we  introduce concrete examples of periodic functions and
 construct their self-intersection and local singular points.
 
In sections 12-14 we introduce and study two  new classes of periodic
functions: the class of periodic helices and the class of S-torus knots.

It is well known  that the folds play an important role in the theory of
the plane curves \cite{Wh}. In the last section,
we construct  folds for the $2$-member chains.
 
 In Appendix A, we prove Proposition \ref{Proposition 11.8} which is essential for Example~\ref{Example 11.1}.
  The last section (Appendix B) contains plots  that correspond to the examples of section~11.
  
  We use notations $F^{\prime}(x)=\big(dF/dx\big)(x)$,  $G^{\prime}(y)=\big(dG/dy\big)(y)$, $\dot{x}=dx/dt$.
 \section{Trigonometric formulas}
  In the present section, we represent $\sin({\ell}x)$ and $\cos({\ell}x)$ $(\ell\in \BN)$ as polynomials of $\sin(x)$ and  $\cos(x)$.  We start with two well known trigonometric formulas:
 \begin{align}&\sin[(\ell+2)x]=\sin(\ell{x})\cos(2x)+\cos(\ell{x})\sin(2x),\quad \ell{\geq}1;\label{2.1}\\ &
 \cos[(\ell+2)x]=\cos(\ell{x})\cos(2x)-\sin(\ell{x})\sin(2x),\quad \ell{\geq}1.\label{2.2}\end{align}
 \begin{Tm}\label{Theorem 2.1}Let  $\ell\in \BN$ be odd, that is, let  $\ell=2n+1$. Then, we have
 \begin{align}&\sin[(2n+1)x]=P_n(\sin^{2}(x))\sin(x),\quad  n{\geq}0;\label{2.3}
 \\ & \cos[(2n+1)x]=(-1)^{n}P_n(\cos^{2}(x))\cos(x),\quad n{\geq}0, \label{2.4}\end{align}
 where $P_n(t)$ is a polynomial of degree $n$.\end{Tm}
 \begin{proof}
Relations \eqref{2.3} and \eqref{2.4} are valid for $n=0.$ It follows from
\eqref{2.1}--\eqref{2.4} that
\begin{equation}P_{n+1}(t^2)=P_n(t^2)(1-2t^2)+2(-1)^{n}P_n(1-t^2)(1-t^2).\label{2.5}\end{equation}
Formulas \eqref{2.1}, \eqref{2.2} and \eqref{2.5} imply
 that 
the assertion of the theorem is valid  for $n+1$ if it is valid for $n$.
\end{proof}
\begin{Rk}\label{Remark 2.2} Using the recurrent formula \eqref{2.5}, one may find explicit expressions for $P_n(t^2)$.\end{Rk}
\begin{Ee}\label{Example 2.3} Clearly, we have
\begin{equation} P_0(t^2)=1,\quad P_1(t^2)=-4t^2+3,\quad P_2(t^2)=16t^4-20t^2+5.\label{2.6}\end{equation}\end{Ee}
\begin{Ee}\label{Example 2.4} In view of \eqref{2.3}, \eqref{2.4}  and \eqref{2.6}, we obtain
\begin{align}& \sin(3x) = -4 \sin^{3}(x) + 3 \sin(x), \quad  \cos(3x) = 4 \cos^{3}(x) - 3 \cos(x),
\label{2.7}
\\ &
 \sin(5x)=16\sin^{5}(x)-20\sin^3(x)+5\sin(x),\label{2.8}
\\ & \cos(5x)=-16\cos^{5}(x)+20\cos^3(x)-5\cos(x).\label{2.9}\end{align}\end{Ee}
\begin{Tm}\label{Theorem 2.5}Let  $\ell\in \BN$ be even, that is, let $\ell=2n$. Then, we have
 \begin{align}&\sin(2nx)=R_n(\sin^{2}(x))\sin(x)\cos(x),\quad n{\geq}1; \label{2.10}
 \\ &
 \cos(2nx)=G_n(\cos^{2}(x)),\quad n{\geq}1, \label{2.11}\end{align}
 where $R_n(t)$ and $G_n(t)$ are polynomials of degrees $n-1$ and $n$  respectively.\end{Tm}
 \begin{proof} Relations \eqref{2.10} and \eqref{2.11} are valid for $n=1.$ It follows from
\eqref{2.1}, \eqref{2.2} and \eqref{2.10}, \eqref{2.11} that
\begin{align}& R_{n+1}(t^2)=R_n(t^2)(1-2t^2)+2G_n(1-t^2),\label{2.12}
\\ &
G_{n+1}(t^2)=-G_n(t^2)(1-2t^2)-2R_n(1-t^2)(1-t^2)t^2.\label{2.13}\end{align}
Formulas \eqref{2.1}, \eqref{2.2} and \eqref{2.12}, \eqref{2.13} imply that 
the assertion of the theorem is valid  for $n+1$ if it is valid for $n$.
\end{proof}
\section{Self-intersection points}
In this section, we study the self-intersection points of the curve  \eqref{1.3}.
\begin{Dn}\label{Definition 3.1} The curve \eqref{1.3} has  a self-intersection at some point $(x_0,y_0)$ if  the curve 
has different tangent lines at this point.\end{Dn}
We rewrite parametric equations \eqref{1.3} in one of the two forms:
\begin{equation} {\mathrm{either}} \quad y=F(x)\quad {\mathrm{or}} \quad x=G(y). \label{3.1}\end{equation}
\begin{Pn}\label{Proposition 3.2}The curve \eqref{1.3} has in the point $(x_0,F(x_0))$ different tangent lines if and only if the derivative $F^{\prime}(x_0)$
has different values.\end{Pn}
\begin{Pn}\label{Proposition 3.3}The curve \eqref{1.3} has in the point $(G(y_0),y_0)$ different tangent lines if and only if the derivative $G^{\prime}(y_0)$
has different values.\end{Pn}
Thus, we need  to calculate the derivatives  $F^{\prime}(x)$ and $G^{\prime}(y)$:
\begin{equation} F^{\prime}(x)={\dot{\psi}(t)}\big/{\dot{\phi}(t)},\quad G^{\prime}(y)={\dot{\phi}(t)}\big/{\dot{\psi}(t)}.\label{3.2}\end{equation}
Using \eqref{2.3}, \eqref{2.4} and \eqref{2.10}, \eqref{2.11} we have
\begin{align}&\phi(t)=\cos(t)Q_{1}(\sin^{2}(t))+U_{1}(\sin^{2}(t)),\label{3.3}
\\ &
\dot{\phi}(t)=\sin(t)[Q_{2}(\sin^{2}(t))+\cos(t)U_{2}(\sin^{2}(t))],\label{3.4}\end{align}
and
\begin{align}&\psi(t)=\sin(t)[S_{1}(\sin^{2}(t))+\cos(t)V_{1}(\sin^{2}(t))],\label{3.5}
\\ &
\dot{\psi}(t)=\cos(t)S_{2}(\sin^{2}(t))+V_{2}(\sin^{2}(t)),\label{3.6}\end{align}
where $Q_k, U_k, S_k, V_k$ ($k=1,2$) are polynomials.
\begin{Rk}\label{Remark 3.4} It follows from \eqref{1.3} that the corresponding curve is symmetric about the axis $y=0$.\end{Rk}

We start by solving the following problem:\\
\emph{Find self-intersection points that belong to the axis $y=0.$}

Assume that  the following conditions are fulfilled for some  $u_0\in \BR$:
\begin{equation}S_1(1-u_{0}^2)+u_0V_1(1-u_{0}^2)=0,\quad 0<u_0<1.\label{3.7}\end{equation}
Now, we introduce  $t_1$ and $t_2$ such that
\begin{equation}  u_0=\cos(t_1)=\cos(t_2),\quad \sin(t_1)=-\sin(t_2)=\sqrt{1-u_{0}^2}>0.\label{3.8}\end{equation}
Using relations \eqref{3.3}--\eqref{3.8}, we obtain the equalities
\begin{align}
& x_0:=\phi(t_1)=\phi(t_2), \quad y_0:=\psi(t_1)=\psi(t_2)=0; \label{3.9-}
\\ &
\dot{\phi}(t_1)=-\dot{\phi}(t_2),\quad \dot{\psi}(t_1)=\dot{\psi}(t_2).\label{3.9}\end{align}
Taking into account  \eqref{3.2} and \eqref{3.9}, we derive the following lemma.
\begin{La}\label{Lemma 3.5} Suppose that 
\begin{equation}\label{3.10n}\dot{\phi}(t_1)\not=0, \quad  \dot{\psi}(t_1)\not=0,  \end{equation}
and that \eqref{3.7} and \eqref{3.8} hold.
Then, the point $(x_0,0)$, where $x_0$ is given in \eqref{3.9-}, is a  self-intersection point of the corresponding curve \eqref{1.3}.\end{La}
\begin{proof} According to \eqref{3.9-},  $t_1$ and $t_2$ generate the same point $(x_0,0)$.
From \eqref{3.2} and \eqref{3.9} we see that
\begin{align}& y^{\prime}(x_0)={\dot{\psi}(t_1)}\big/{\dot{\phi}(t_1)} \quad {\mathrm{for}} \quad t=t_1,\label{3.11}
\\ &
y^{\prime}(x_0)=-{\dot{\psi}(t_1)}\big/{\dot{\phi}(t_1)} \quad {\mathrm{for}} \quad t=t_2.\label{3.12}\end{align}
Thus, the curve has two different tangent lines.\end{proof}
Our remarks below deal with the cases where \eqref{3.10n} is not valid.
\begin{Rk}\label{Remark 3.6} If $\dot{\phi}(t_1)=0$  but
$\dot{\psi}(t_1)\not=0$, 
then, according to \eqref{3.11} and \eqref{3.12}, we have $y^{\prime}(x_0)=\infty$
at the both points $t_1$ and $t_2$. Hence, the point $(x_0,0)$ is not   a self-intersection point of the corresponding curve \eqref{1.3}.
\end{Rk}

\begin{Rk}\label{Remark 3.7}
 If $\dot{\psi}(t_1)=0$  but
$\dot{\phi}(t_1)\not=0$, 
then, according to \eqref{3.11} and \eqref{3.12}, we have $y^{\prime}(x_0)=0$
at the both points $t_1$ and $t_2$. Hence, the point $(x_0,0)$ is not   a self-intersection point of the corresponding curve \eqref{1.3}.
\end{Rk}

\begin{Rk}\label{Remark 3.8}If $\dot{\phi}(t_1)=0$ and $\dot{\psi}(t_1)=0$, 
then
the point $(x_0,0)$ is a local singular point.
\end{Rk}
Lemma \ref{Lemma 3.5} and Remarks \ref{Remark 3.6}--\ref{Remark 3.8} imply the next theorem.
\begin{Tm}\label{Theorem 3.9} Suppose  that the relations \eqref{3.7}, where $u_0=\cos(t_1)$, and \eqref{3.8} hold. Set $x_0:=\phi(t_1)$.
Then, the point $(x_0,0)$ is a   self-intersection point of the corresponding curve \eqref{1.3} if
and only if  the inequalities \eqref{3.10n} are fulfilled.\end{Tm}

Next, we consider the  self-intersection problem for another axis:\\
\emph{Find self-intersection points that belong to the axis $x=0.$}

We study the case where all numbers $m_k$ are odd.
Using relations \eqref{2.3} and \eqref{2.4}, we (similar to \eqref{3.3}--\eqref{3.6})  obtain the equalities:
\begin{align}& \phi(t)=\cos(t)Q_{1}(\sin^{2}(t)),\quad \dot{\phi}(t)=\sin(t)Q_{2}(\sin^{2}(t));\label{3.18}
\\ &
\psi(t)=\sin(t)S_{1}(\sin^{2}(t)),\quad \dot{\psi}(t)=\cos(t)S_{2}(\sin^{2}(t)).\label{3.19}\end{align}

We assume that for some real $u_0$ the following conditions are fulfilled:
\begin{equation}Q_1(1-u_0^2)=0,\quad 0<u_0<1,\label{3.20}\end{equation}
and introduce  such $t_1$ and $t_2$ that
\begin{equation}  u_0=\cos(t_1)=-\cos(t_2),\quad \sin(t_1)=\sin(t_2)=\sqrt{1-u_0^2}\, .\label{3.21}\end{equation}
Using relations \eqref{3.18}--\eqref{3.21}, we obtain the equalities
\begin{align}&
 x_0:=\phi(t_1)=\phi(t_2)=0, \quad y_0:=\psi(t_1)=\psi(t_2); \label{3.22-}
\\ &
\dot{\phi}(t_1)=\dot{\phi}(t_2),\quad \dot{\psi}(t_1)=-\dot{\psi}(t_2).\label{3.22}\end{align}
Taking into account relations \eqref{3.2} and \eqref{3.22} we get:
\begin{La}\label{Lemma 3.10} Suppose that 
\begin{equation}\label{3.23}\dot{\phi}(t_1)\not=0, \quad  \dot{\psi}(t_1)\not=0,  \end{equation}
and that the relations \eqref{3.20} and \eqref{3.21} hold.
Then, the point $(0,y_0)$, where $y_0$ is given in \eqref{3.22-}, is a  self-intersection point of the corresponding curve \eqref{1.3}.\end{La}
\begin{proof} According to \eqref{3.22-},  $t_1$ and $t_2$ generate the same point $(0,y_0)$.
From \eqref{3.2} and \eqref{3.22} we see that
\begin{align}& x^{\prime}(y_0)={\dot{\phi}(t_1)}\big/{\dot{\psi}(t_1)} \quad {\mathrm{for}} \quad t=t_1,\label{3.24}
\\ &
x^{\prime}(y_0)=-{\dot{\phi}(t_1)}\big/{\dot{\psi}(t_1)} \quad {\mathrm{for}} \quad t=t_2.\label{3.25}\end{align}
Thus, the curve has two different tangent lines.\end{proof}
Similar to the case of Lemma \ref{Lemma 3.5} and Theorem \ref{Theorem 3.9}, Lemma 3.10 implies the following theorem
\begin{Tm}\label{Theorem 3.11} Suppose  that the relations \eqref{3.20}, where $u_0=\cos(t_1)$, and \eqref{3.21} hold. 
Set $y_0:=\psi(t_1)$. Then, the point $(0,y_0)$ is a    self-intersection point of the corresponding curve \eqref{1.3} if
and only if  the inequalities \eqref{3.23} are fulfilled. \end{Tm}
\begin{Pn}\label{Proposition 3.12}If all  $m_k$ are odd, then the corresponding curve \eqref{1.3}
is symmetric about the axis $x=0$.\end{Pn}
\begin{Rk}\label{RkS} The self-intersection points on the axes $y = 0$ or $x = 0$ often appear as a result of the symmetry of the corresponding
curve.
\end{Rk}
\section{$N$-member chains: the subclass \eqref{1.1}}
\emph{In this section, we  consider  the special case \eqref{1.1}  of the $n$-member chains.}

Let us represent a curve $f(t)$ of the   subclass \eqref{1.1} in the form:
\begin{equation} f(t)=e^{itm_1}\Big[c_1+\sum_{k=2}^{n}c_{k}e^{it(m_k-m_1)}\Big],\label{4.1}\end{equation}
where $m_1$  is the smallest among the numbers $m_k$. 
Hence, we have the following periodicity of $|f(t)|$ and $\arg\big(f(t)\big)$:
\begin{align}& r(t):=|f(t)|=r(t+2\pi/Q),\label{4.2}
\\ &
\alpha(t):=\arg\Big[c_1+\sum_{k=2}^{n}c_{k}e^{it(m_k-m_1)}\Big]=\alpha(t+2\pi/Q),\label{4.3}\end{align}
where $Q$ is the greatest common divisor of the numbers $m_k-m_1$  $(2{\leq}k{\leq}n)$.
\begin{Pn}\label{Proposition 4.1}If the point $w_0=f(t_0)$ is a self-intersection point of the curve $f(t)$, then the points
\begin{equation} w_k:=f(t_0+2k\pi/Q)\quad (0{\leq}k{\leq}Q-1) \label{4.4}\end{equation}
are self-intersection points of the curve $f(t)$ as well. \end{Pn}
\begin{proof} The conditions of the proposition imply that we have 
\begin{align}& \label{z0}
f(t_0)=f(s_0), \quad \dot{y}(t_0)\big/ \dot{x}(t_0)\not= \dot{y}(s_0)\big/ \dot{x}(s_0) 
\end{align}
for $x(t):=\Re\big(f(t)\big), \,\, y(t):=\Im\big(f(t)\big)$ and for some $s_0\not=t_0$. In view of the equality in \eqref{z0}, relations \eqref{4.1}--\eqref{4.3} yield
\begin{equation}f(t_0+2k\pi/Q)=e^{im_1(2k\pi/Q)}f(t_0)   =e^{im_1(2k\pi/Q)}f(s_0)=f(s_0+2k\pi/Q).\label{4.6}\end{equation}
Moreover, it follows from  \eqref{4.1} that $\dot{f}\Big(t+\frac{2k\pi}{Q}\Big)=\exp\{2 i m_1 k \pi \big/Q\}\dot{f}(t)$.
Hence,  we obtain:
\begin{align}& \label{z1}\dot{x}(t_0+2k\pi/Q)=\cos(\b)\dot{x}(t_0)-\sin(\b)\dot{y}(t_0), \quad \b:=m_1(2k\pi/Q);
\\ & \label{z2}
\dot{y}(t_0+2k\pi/Q)=\sin(\b)\dot{x}(t_0)+\cos(\b)\dot{y}(t_0).\end{align}
From \eqref{z1} and \eqref{z2}, we derive
\begin{align}& \label{z3}
\dot{y}(t_0+2k\pi/Q)\big/ \dot{x}(t_0+2k\pi/Q)=\frac{\sin(\b)+\cos(\b)g_1}{\cos(\b)-\sin(\b)g_1},  \quad g_1:=\dot{y}(t_0)\big/ \dot{x}(t_0).\end{align}
In a similar way, we derive
\begin{align}& \label{z4}
\dot{y}(s_0+2k\pi/Q)\big/ \dot{x}(s_0+2k\pi/Q)=\frac{\sin(\b)+\cos(\b)g_2}{\cos(\b)-\sin(\b)g_2},  \quad g_2:=\dot{y}(s_0)\big/ \dot{x}(s_0).\end{align}
Since $g_1\not= g_2$, relations \eqref{z3} and \eqref{z4} yield
\begin{align}& \label{z5}
\dot{y}(t_0+2k\pi/Q)\big/ \dot{x}(t_0+2k\pi/Q)\not= \dot{y}(s_0+2k\pi/Q)\big/ \dot{x}(s_0+2k\pi/Q).\end{align}
 \end{proof}
\section{$2$-member chains: the subclass \eqref{1.1}}
In this section, we consider the $2$-member chains:
\begin{equation} f(t)={c_1}e^{itm}+{c_2}e^{it\ell}=e^{itm}[c_1+c_2e^{it(\ell-m)}],\label{5.1}\end{equation}
where $m$ and $\ell$ are coprime integers and $\ell>m.$ Clearly, the greatest common divisor of  $\ell-m$ is $Q=\ell-m.$
Thus, relations \eqref{4.2}  and \eqref{4.3} take the form:
\begin{align}& r(t):=|f(t)|=r(t+2\pi/(\ell-m)),\label{5.2}
\\ &
\alpha(t):=\arg[c_1+c_2e^{it(\ell-m)}]=\alpha(t+2\pi/(\ell-m)).\label{5.3}\end{align}
Proposition \ref{Proposition 4.1} implies:
\begin{Pn}\label{Proposition 5.1}If the point $w_0=f(t_0)$ is a self-intersection point of the curve $f(t)$, then the points
\begin{equation} w_k=f(t_0+2k\pi/(\ell-m))\quad (0{\leq}k{\leq}\ell-m-1) \label{5.4}\end{equation}
are self-intersection points of the curve $f(t)$ as well. \end{Pn}
Let us compare the curves
\begin{equation}f(t)={c_1}e^{itm}+{c_2}e^{it\ell}\quad {\mathrm{and}}\quad g(t)={c_1}e^{itm}-{c_2}e^{it\ell}.\label{5.5}\end{equation}
Using formula \eqref{5.1}, we derive
\begin{equation}f\Big(t+\frac{\pi}{\ell-m}\Big)=e^{im\pi/(\ell-m)}g(t).\label{5.6}\end{equation}
The next assertion follows directly from \eqref{5.5} and \eqref{5.6}.
\begin{Cy}\label{Corollary 5.2} Rotating the curve $g(t)$ by an angle $m\pi/(\ell-m)$ around the origin,  we obtain
the curve $f\Big(t+\frac{\pi}{\ell-m}\Big)$.\end{Cy}
Formula \eqref{5.1} implies:
\begin{equation}r^{2}(t)=c_{1}^2+c_{2}^2+2c_{1}c_{2}\cos[(\ell-m)t],\label{5.7}\end{equation}
where $r(t)=|f(t)|$ as in \eqref{5.2}.
\begin{Pn}\label{Proposition 5.3} Let $r(t_1)=r(t_2)$. Then, either
\begin{equation} t_1=t_2+2s\pi/(\ell-m) \quad ({\mathrm{mod}} \,\, 2\pi),\quad 0{\leq}s{\leq}\ell-m-1, \label{5.8}\end{equation}
or
\begin{equation} t_1=-t_2+2s\pi/(\ell-m)\quad ({\mathrm{mod}} \,\, 2\pi),\quad 0{\leq}s{\leq}\ell-m-1. \label{5.9}\end{equation}\end{Pn}
\begin{proof}
It follows from \eqref{5.7} that
\begin{equation}\cos[(\ell-m)t_1]-\cos[(\ell-m)t_2]=0.\label{5.10}\end{equation}
Hence, we have
\begin{equation}2\sin[(t_1+t_2)(\ell-m)/2]\sin[(t_1-t_2)(\ell-m)/2]=0.\label{5.11}\end{equation}
The statement of proposition directly  follows  from \eqref{5.11}.
\end{proof}
Let us consider self-intersections at the points $\om_0\not=0$ (the case $\om_0=0$ will be treated separately).
\begin{Tm}\label{Theorem 5.4}
If 
$$w_0=f(t_0)=f(s_0)\not=0\quad \big(t_0\not=s_0 \,\, \mod(2\pi)\big)$$ is a self-intersection point of the curve \eqref{5.1},
then the following assertions are valid:\\
1. For $$t_k:=t_0+2k\pi/(\ell-m)\,\, {\mathrm{and}}\,\,  s_j:=s_0+2j\pi/(\ell-m)\quad (0{\leq}k, j{\leq}\ell-m-1),$$
we have
\begin{align}& \label{5.14}
t_k\not=t_j  \mod(2\pi), \quad s_k\not=s_j  \mod(2\pi) \quad {\mathrm{for}} \quad k\not=j;
\\ & \label{5.14+}
t_k\not=s_j  \mod(2\pi) \quad {\mathrm{for\,\, all}} \quad   0{\leq}k, j{\leq}\ell-m-1.
\end{align}
2. The following equalities are fulfilled:
\begin{equation}f(t_k)=f(s_k),\quad 0{\leq}k{\leq}\ell-m-1. \label{5.13}\end{equation}
3. The points $w_k$ of the form
\begin{equation}w_k=f(t_k)=f(s_k)\quad (0{\leq}k{\leq}\ell-m-1) \label{5.12}\end{equation}
are   different self-intersection points of the curve \eqref{5.1}.\\
4. If a point $w=f(\tau)$ of the curve \eqref{5.1} satisfies the condition $|w|=|f(t_0)|$ then we have
\begin{equation} \tau=t_k \quad ({\mathrm{mod}} \,\, 2\pi)\quad {\mathrm{or}} \quad \tau= s_j  \mod(2\pi) \quad  (0{\leq}k, j{\leq}\ell-m-1).\label{5.14!}\end{equation}
5. If $\tau$ satisfies \eqref{5.14!}, then  $(-\tau)$ satisfies \eqref{5.14!} as well.\\
6.  For the self-intersection point $\om_k$ $($see \eqref{5.12}$)$ there is only one point, namely, point $s_k \mod(2\pi)$, such that $f(t_k)=f(s_k)$ and $s_k\not= t_k \mod(2\pi)$.
 \end{Tm}
\begin{proof} Cleary \eqref{5.14} holds. Since $t_0\not=s_0\mod(2\pi)$ but $f(t_0)=f(s_0)$, it follows that   $t_k\not=s_k\mod(2\pi)$ but $f(t_k)=f(s_k)$ (see \eqref{5.1}).
Thus, 
\begin{equation}f(t_k)\not=f(s_j) \,\,  {\mathrm{for}}\,\, k\not=j, \,\, {\mathrm{although}} \,\, |f(t_k)|=|f(s_j)|\label{5.15-}\end{equation}
(see again \eqref{5.1}).
It proves \eqref{5.14+} for both cases $k=j$ and $k\not=j$. As mentioned above, the equalities \eqref{5.13}
(assertion 2) are valid as well.

Assertion 3 follows from Proposition \ref{Proposition 5.1} and assertion 2.

According to Proposition \ref{Proposition 5.3}, the equation $|w_0|=r(\tau)$ has no more than $2(\ell - m)$ solutions (up to $\mod(2\pi)$).
On the other hand, in view of  \eqref{5.1} and \eqref{5.15-}, we have $|f(t_k)|=|f(t_0)|=|f(s_j)$.
Thus, relations \eqref{5.14!} present solutions of $|w_0|=r(\tau)$. Taking into account  assertion 1, we see that
all $2(\ell-m)$ solutions are presented in \eqref{5.14!}.

It follows from \eqref{5.1} that $f(-\tau)=\ov{f(\tau)}$. Hence, the equalities in \eqref{5.14!} and \eqref{5.15-} yield $|f(-\tau)|=|f(\tau)|=|f(t_0)|$.
Therefore, the proof of assertion 5 follows from the proof of assertion 4.

Assertion 4 and the inequality in \eqref{5.15-} yield assertion 6.
\end{proof}
Next, we consider the case $w_0=f(t_0)=0$.
It is easy to see that the proposition below is valid.
\begin{Pn}\label{Proposition 5.5} Let $f(t)$ have the form \eqref{5.1}. If $f(t_0)=0$, then either $c_1=c_2$ or $c_1=-c_2.$\end{Pn}
We note that the equality $f(t_0)=0$  is equivalent to the equality $r(t_0)=0$. Therefore, taking into account relation \eqref{5.7} 
we obtain the following result for the cases $c_1=c_2$ and $c_1=-c_2.$
\begin{Tm}\label{Theorem 5.6} Let $c_1=-c_2$. Then, $f(t_k)=0$ if and only if the point $t_k$ has the form
\begin{equation}t_k=2k\pi/(\ell-m)\mod(2\pi), \quad  0{\leq}k{\leq}\ell-m-1. \label{5.16}\end{equation}
 Let $c_1=c_2.$ Then $f(t_k)=0$ if and only if the point $t_k$ has the form
\begin{equation}t_k=(2k+1)\pi/(\ell-m) \mod(2\pi), \quad  0{\leq}k{\leq}\ell-m-1. \label{5.17}\end{equation}\end{Tm}
\begin{Dn}\label{Definition 5.7} We say that the self-intersection points $w_k=f(t_k)$ belong to the same module class if relations \eqref{5.13}
are fulfilled.\end{Dn}
\begin{Cy}\label{Corollary 5.8} Let the conditions of  Theorem \ref{Theorem 5.4} be fulfilled.
If the number $(\ell-m)$ is odd, then one and only one member of the module class \eqref{5.12} is real.\end{Cy}
\begin{proof} From Theorem \ref{Theorem 5.4} and Remark \ref{Remark 3.4} we have the assertions:\\
1. The number of points $w_k$ is equal to $(\ell-m)$ .  Hence, this number is odd.\\
2. The number of the points $w_k$ outside the axis $y=0$ is even.\\
3. The number of the points $w_k$ which belong to the axis $y=0$ cannot be greater then two.\\
The statement of the corollary follows from these assertions.
\end{proof}
\begin{Cy}\label{Corollary 5.9} Let the conditions of  Theorem \ref{Theorem 5.4} be fulfilled.
If $(\ell-m)=2s$ and $s$ is odd, then either two and only two members of module class \eqref{5.12} belong to the
axis $y=0$ or  two and only two members of module class \eqref{5.12} belong to the
axis $x=0$.\end{Cy}
\begin{proof} The curve $w=f(t)$  is symmetric about the axes $y=0$ and $x=0$.Therefore, the number of members of the module class \eqref{5.12}
that do not belong to the axis $y=0$ and $x=0$ is a multiple of 4. On each of the axes $y=0$ and $x=0$ two and only two points can belong  to the  module class \eqref{5.12}.
\end{proof}
\begin{Rk}\label{Remark 5.10} If we know one point of the module class \eqref{5.12}, we know all the points of this module class.
A method for constructing self-intersection points
belonging either $y=0$ or $x=0$  is given in section 3 of the present paper.\end{Rk}

The next assertion follows from \eqref{5.1} and \eqref{5.12}.
\begin{Cy}\label{Corollary 5.11} Let the conditions of Theorem \ref{Theorem 5.4} be fulfilled. Then,
\begin{equation}w_k=w_0e^{2 i k\pi{m}/(\ell-m)} \mod(2\pi),\quad 0{\leq}k{\leq}\ell-m-1.\label{5.18}\end{equation}
 \end{Cy}
 Let us consider the finite group of rotations of the complex plane with multiplications by $g_k$ as elements, where
\begin{equation}g_k=e^{2i k\pi{m}/(\ell-m)},\quad 0{\leq}k{\leq}\ell-m-1.\label{5.19}\end{equation} 
We denote this group by $G(\ell,m).$ Relations \eqref{5.18} and \eqref{5.19} yield the next theorem.
\begin{Tm}\label{Theorem 5.12}  Let the conditions of  Theorem \ref{Theorem 5.4} be fulfilled. Then, the elements of the group $G(\ell,m)$ map the corresponding 
module classes of self-intersection points onto themselves.\end{Tm}
\begin{Rk}\label{Remark 5.13}It is well-known that the  groups of rotations of the complex plane
are Lie groups. In particular,  $G(\ell,m)$ is a Lie group.\end{Rk}

\section{Local singular points}In this section, we investigate the local singular points of the curve \eqref{1.3}. We use the following assertion (see \cite{BrG}).
\begin{Pn}\label{Proposition 6.1}
The local singular points of the curve \eqref{1.3} are those points where 
\begin{equation}\frac{d\phi(t)}{dt}=\frac{d\psi(t)}{dt}=0.\label{6.1}\end{equation}
\end{Pn}
It follows from \eqref{1.3} and \eqref{6.1} that
\begin{align}&\frac{d\phi(t)}{dt}=-\sum_{k=1}^{n}m_{k}c_{k}\sin(m_{k}t)=0.\label{6.2}
\\ &
\frac{d\psi(t)}{dt}=-\sum_{k=1}^{n}m_{k}d_{k}\cos(m_{k}t)=0.\label{6.3}\end{align}
Taking into account relations \eqref{2.3}, \eqref{2.4}, \eqref{2.10}, \eqref{2.11} as  well as the equalities \eqref{6.2} and \eqref{6.3} above, we have
\begin{align}&\frac{d\phi(t)}{dt}=A(\sin(t))\cos(t)+B(\sin(t))=0,\label{6.4}
\\ &
\frac{d\psi(t)}{dt}=C(\sin(t))\cos(t)+D(\sin(t))=0,\label{6.5}\end{align}
where $A(x),\, B(x),\, C(x)$ and $D(x)$ are polynomials. Relations \eqref{6.4} and \eqref{6.5} 
written down in the matrix form $\begin{bmatrix}A&B\\C&D\end{bmatrix}\begin{bmatrix}\cos(t)\\1 \end{bmatrix}=0$
imply that:
\begin{equation}A(\sin(t))D(\sin(t))-B(\sin(t))C(\sin(t))=0.\label{6.6}\end{equation}
Assuming that relations
\begin{equation}A(u_1)D(u_1)-B(u_1)C(u_1)=0,\quad 0{\leq}u_1{\leq}1\label{6.7}\end{equation}
hold for some $u_1$ and fixing  $t_1$ such that
\begin{equation} u_1=\sin(t_1),\label{6.8}\end{equation}
we  formulate the main result of this section.
\begin{Tm}\label{Theorem 6.2} Let $x_1=\phi(t_1)$ and $y_1=\psi(t_1)$. Then, the point $(x_1,y_1)$ is
a local singular point of the curve \eqref{1.3} if and only if  the relations \eqref{6.7}, \eqref{6.8} and the
equality \eqref{6.4} $($for $t_1):$
\begin{equation}A(\sin(t_1))\cos(t_1)+B(\sin(t_1))=0\label{6.9}\end{equation}
are fulfilled.\end{Tm}
\section{Local singular points: the case $n=2$}
1. Let us consider the curve \eqref{1.1} with $n=2$. That is, we assume that
\begin{equation} f(t)={c_1}e^{itm}+{c_2}e^{it\ell},\quad  \ell>m,\quad \ell\not=-m,\label{7.1}\end{equation}
where $m$ and $\ell$ are coprime integers.
In the singular points of the curve \eqref{7.1}, we have
\begin{equation}\dot{\phi}(t)=-mc_{1}\sin(tm)-{\ell}c_{2}\sin(t\ell)=0,\label{7.2}\end{equation}
\begin{equation}\dot{\psi}(t)=mc_{1}\cos(tm)+{\ell}c_{2}\cos(t\ell)=0.\label{7.3}\end{equation}
Relations \eqref{7.2} and \eqref{7.3} imply that $(mc_1)^2\sin^2(tm)=(\ell{c_2})^2\sin^2(t\ell)$ and
$(mc_1)^2\cos^2(tm)=(\ell{c_2})^2\cos^2(t\ell)$. Therefore, we obtain the next proposition.
\begin{Pn}\label{Proposition 7.1} If the curve \eqref{7.1} has a local singular point, then 
\begin{equation} (mc_1)^2=(\ell{c_2})^2.\label{7.4}\end{equation}
\end{Pn}
We introduce double singularities, that is, singular points where
\begin{equation}[\ddot{\phi}(t)]^2+[\ddot{\psi}(t)]^2\not=0.\label{7.5}\end{equation}
\begin{Pn}\label{Proposition 7.2} All the local singular points of the  curve \eqref{7.1} are double. 
\end{Pn}
\begin{proof} Suppose that $\ddot{\phi}(t)=\ddot{\psi}(t)=0$ in a local singular point $t$. Then, we have
\begin{align}& m^2{c_1}\cos(tm)+\ell^{2}{c_2}\cos(t\ell)=0, \label{7.6}
\\ &
m^2{c_1}\sin(tm)+\ell^{2}{c_2}\sin(t\ell)=0. \label{7.7}\end{align}
Using \eqref{7.6} and \eqref{7.7}, we (similar to \eqref{7.4}) derive  that $(m^2{c_1})^2=(\ell^{2}{c_2})^2$.
Thus, taking into account \eqref{7.4}, we obtain
$\ell^2=m^2,$ which contradicts the condition \eqref{7.1}. 
\end{proof}

2. Now,  consider separately the case where
\begin{equation}m{c_1}=\ell{c_2}.\label{7.8}\end{equation}
Taking into account \eqref{7.2}, \eqref{7.3} and \eqref{7.8},  we have
\begin{align}& \nn \dot{\phi}(t)=-mc_{1}[\sin(tm)+\sin(t\ell)]=0,\\ & \label{7.9}\\
& \nn
 \dot{\psi}(t)=mc_{1}[\cos(tm)+\cos(t\ell)]=0.\end{align}
Equalities \eqref{7.9} may be rewritten in the form
\begin{align}&\nn \dot{\phi}(t)=-2m{c_1}\sin\left(\frac{\ell+m}{2}t\right)\cos\left(\frac{\ell-m}{2}t\right)=0,  
\\ \label{7.10} 
\\ & \nn
\dot{\psi}(t)=2m{c_1}\cos\left(\frac{\ell+m}{2}t\right)\cos\left(\frac{\ell-m}{2}t\right)=0. \end{align}
Relations \eqref{7.10} are equivalent to the equality
\begin{equation} \cos\left(\frac{\ell-m}{2}t\right)=0. \label{7.11} \end{equation}
Therefore, the relations \eqref{7.10} are valid if and only if
\begin{equation}t_k=\frac{2k+1}{\ell-m}\pi \,\mod(2\pi),\quad 0{\leq}k{\leq}\ell-m-1.\label{7.12}\end{equation}
\begin{Tm}\label{Theorem 7.3} Let the condition \eqref{7.8} be fulfilled. Then,
 the set of the local singular points of the curve \eqref{7.1} coincides with the set  $w_k=f(t_k),$
 where $t_k$ are given in \eqref{7.12}. \end{Tm}
3. Next, consider the case
\begin{equation}m{c_1}=-\ell{c_2}.\label{7.13}\end{equation}
Taking into account \eqref{7.2}, \eqref{7.3} and \eqref{7.13}  we have
\begin{align}&\dot{\phi}(t)=-mc_{1}[\sin(tm)-\sin(t\ell)]=0,\nonumber
\\ & \label{7.14}
\\ &\nn
\dot{\psi}(t)=mc_{1}[\cos(tm)-\cos(t\ell)]=0.\end{align}
Equalities \eqref{7.14} can be written in the form
\begin{align}& \dot{\phi}(t)=2m{c_1}\sin\left(\frac{\ell-m}{2}t\right)\cos\left(\frac{\ell+m}{2}t\right)=0, \nonumber
\\ &  \label{7.15}
\\ \nn &
\dot{\psi}(t)=2m{c_1}\sin\left(\frac{\ell+m}{2}t\right)\sin\left(\frac{\ell-m}{2}t\right)=0. \end{align}
Relations \eqref{7.15} are equivalent to the equality
\begin{equation} \sin\left(\frac{\ell-m}{2}t\right)=0. \label{7.16} \end{equation}
Therefore, the relations \eqref{7.15} are valid if and only if
\begin{equation}t=\tau_k=\frac{2k}{\ell-m}\pi \mod(2\pi),\quad 0{\leq}k{\leq}\ell-m-1.\label{7.17}\end{equation}
Hence, we proved the following  assertion.
\begin{Tm}\label{Theorem 7.4} Let the condition \eqref{7.13} be fulfilled. Then,
the set of the local singular points of the curve \eqref{7.1} coincides with the set  $w_k=f(\tau_k).$ \end{Tm}
4. Let us introduce the notion of the return point.
\begin{Dn}\label{Definition 7.5} The local singular point $t_0$ is a return point of the first kind if the following conditions are fulfilled.\\
I. Both branches $(t>t_0$ and $t<t_0)$ of the curve $w=f(t)$ have a common tangent line at the singular point $t_0$.\\
II. For small $(t-t_0)$  both branches of the curve $w=f(t)$ are located on a one side of the common normal.\\
III. For small $(t-t_0)$  both branches of the curve $w=f(t)$ are located on different sides of the common tangent line.\end{Dn}
\begin{Tm}\label{Theorem 7.6} All the local singular points of the curve \eqref{7.1} are return points of the first kind.\end{Tm}
\begin{proof} According to Proposition \ref{Proposition 7.1}, the existence of the local singular points yields \eqref{7.4}
First, consider the case \eqref{7.13}. We assume that $k=0$. In view of \eqref{7.14} and \eqref{7.17} we have:
\begin{align}& \tau_0=0,\quad x(t)=(c_1+c_2)-t^{2}c_{1}m(m+\ell)/2 +O(t^4),\label{7.18}
\\ &
 y(t)=-t^{3}c_{1}m(m^2+\ell^2)/6 +O(t^5).\label{7.19}
 \end{align}
It follows from \eqref{7.13} and the inequality $m\not=\pm\ell$ that
\begin{equation}c_1+c_2\not=0,\quad m+\ell\not=0,\quad m^2+\ell^2\not=0.\label{7.20}\end{equation}
We need also the formula
\begin{equation} y^{\prime}(x)=\tan\left(\frac{\ell+m}{2}t\right), \label{7.21}\end{equation} which follows from \eqref{7.15}.
Using equalities $\tau_0=0$ and \eqref{7.21} we obtain
\begin{equation} y^{\prime}(x_0)=0 \quad (x_0=x(0)). \label{7.22}\end{equation}
Thus,  property I of Definition 7.5 is fulfilled (see Proposition \ref{Proposition 3.2}). Property II follows from \eqref{7.18} . The equality \eqref{7.19}
implies property III. Hence, the local singular point  $(x_0,0)$ of the curve \eqref{7.1} is a return point of the first kind.
All other  local singular point   of the same module class are obtained from the point  $(x_0,0)$  by rotating
 the  the curve \eqref{7.1} by a certain angle
(see \eqref{5.2}, \eqref{5.3}).   In this way our  theorem is proved for the case \eqref{7.13}.

Now, formula \eqref{5.6}, where $f$ and $g$ have the forms \eqref{5.5}, implies that our theorem is valid in the case \eqref{7.8} as well. 
\end{proof}
The next assertion follows from \eqref{7.12} and \eqref{7.17} (for the cases \eqref{7.8} and \eqref{7.13}, respectively).
\begin{Cy}\label{Corollary 7.7} Let the condition \eqref{7.4} be fulfilled. 

Then, for the local singular points we have
\begin{equation}w_k=w_0e^{2 i k\pi/(\ell-m)}\quad (0{\leq}k{\leq}\ell-m-1),\label{7.23}\end{equation}
where $w_0=e^{i\pi/(\ell-m)}$ if the condition \eqref{7.8} is fulfilled  and $w_0=1$  if the condition \eqref{7.13} is fulfilled.
 \end{Cy}
 Recall (see Remark \ref{Remark 5.13}) that  the finite group of  rotations of the complex plane  obtained via the multiplication by
\begin{equation}g_k=e^{2 i k\pi/(\ell-m)}\quad (0{\leq}k{\leq}\ell-m-1)\label{7.24}\end{equation}
is a Lie group denoted by $G(\ell,m)$. (It is easy to see that the families $g_k$ given by \eqref{5.19} and \eqref{7.24} coincide.)
From  \eqref{7.23} and \eqref{7.24} we obtain an assertion.
\begin{Tm}\label{Theorem 7.8} Let the condition \eqref{7.4} be fulfilled. Then, the elements of the group $G(\ell,m)$ map the corresponding
 classes of singular points onto themselves.\end{Tm}

\section{Classical curves}
In this section, we consider the well known curves which belong to the $\break 2$-member chain class.

I. \emph{Epicycloid}\cite{Nas, Pro}.
\begin{Dn}\label{Definition 8.1} An epicycloid is a plane curve traced out by a point on the circumference
 of a circle which rolls without slipping on a fixed circle in the same plain. $($Note that the rolling circle rolls outside of the fixed circle.$)$\end{Dn}
 If the fixed circle has the radius $R$ and the rolling  circle has the radius $r$, then the parametric equations  for the epicycloid have the form:
\begin{align} & x=\phi(t)=(R+r)\cos(t)-r\cos\left(\frac{R+r}{r}t\right),\label{8.1}
\\ &
 y=\psi(t)=(R+r)\sin(t)-r\sin\left(\frac{R+r}{r}t\right).\label{8.2}\end{align}
The epicycloid \eqref{8.1}, \eqref{8.2} may be written in the form
\begin{equation}w=f(t)=(R+r)e^{it\lambda_1}-re^{it\lambda_2},\label{8.4+}\end{equation}
where $\lambda_1=1,\quad \lambda_2=(R+r)/r.$
\begin{Rk}\label{Remark 8.3} If $(R+r)/{r}$ is an irrational  number, then the numbers $\lambda_1$
and $\lambda_2$ are linear independent over the field of rational numbers. Hence, the function $w=f(t)$
is almost periodic. For that case, the curve $w=f(t)$ was  studied in our paper \cite{Sak2}.\end{Rk}

 If $(R+r)/r$ is a rational number, it may be represented in the form\\
 $(R+r)/r=p/q,$ where $p$ and $q$ are positive, coprime integers, $p>q.$ In this case, 
the results  for the curve \eqref{5.1}: 
\begin{equation} f(s)={c_1}e^{ism}+{c_2}e^{is\ell}\label{5.1'}\end{equation}
may be applied. Compare \eqref{8.4+} and \eqref{5.1'} to see that the following proposition is valid.
\begin{Pn}\label{Pn8.2} System \eqref{8.1}, \eqref{8.2},  where
\begin{equation} c_1=R+r,\quad c_2=-r,\quad m=q,\quad \ell=p, \quad t=qs,\label{8.3}\end{equation}
is equivalent to equation \eqref{5.1'}.
\end{Pn}
The self-intersection points of the epicycloid, where  $(R+r)/r$ is rational number, are described
in Theorem \ref{Theorem 5.6} and in Corollary \ref{Corollary 5.8}.  Moreover,
it follows from the relations \eqref{8.3} that $mc_1=-\ell{c_2}$, that is, condition \eqref{7.13} is fulfilled.
Thus, the local singular points of the epicycloid
are described in Theorem \ref{Theorem 7.4}.

II. \emph{Hypocycloid}\cite{Nas, Pro}.
\begin{Dn}\label{Definition 8.4} A hypocycloid is a plane curve traced  by a point on the circumference
 of a circle which rolls without slipping on a fixed circle in the same plain, where the rolling circle rolls inside of the fixed circle.\end{Dn}
If the fixed circle has the radius $R$ and the rolling  circle has the radius $r$, the parametric equations
for the hypocycloid have the form \cite{BrG}:
\begin{align}& x=\phi(t)=(R-r)\cos(t)+r\cos\left(\frac{R-r}{r}t\right),\label{8.6}
\\ &
y=\psi(t)=(R-r)\sin(t)-r\sin\left(\frac{R-r}{r}t\right).\label{8.7}\end{align}
The hypocycloid \eqref{8.6}, \eqref{8.7} may be written in the form
\begin{equation}w=f(t)=(R-r)e^{it\lambda_1}+re^{it\lambda_2},\label{8.10}\end{equation}
where $\lambda_1=1,\quad \lambda_2=-(R-r)/r.$
If $(R-r)/r$ is an irrational number, we deal with the case of almost periodic functions discussed in 
Remark~\ref{Remark 8.3}. 
 If $(R-r)/r$ is a rational number, it may be represented in the form
 $\break (R-r)/r=p/q,$ where $p$ and $q$ are positive, coprime integers. 
\begin{Pn}\label{8.5} System \eqref{8.5}, \eqref{8.6} is equivalent to equation \eqref{5.1'}, where
\begin{equation}  c_1=R-r,\quad  c_2=r\quad m=-q,\quad \ell=p, \quad t=-qs.\label{8.8}\end{equation}\end{Pn}
The self-intersection points of the hypocycloid, where  $(R-r)/r$ is a rational number, are described
in Theorem \ref{Theorem 5.6} and in Corollary \ref{Corollary 5.8}.  
Moreover, it follows from \eqref{8.8} that the equality \eqref{7.13} is valid again. Thus,
the local singular points of the hypocycloid
are described in Theorem \ref{Theorem 7.4}.


III. \emph{Epitrochoid} \cite{Lud, Nas}.
\begin{Dn}\label{Definition 8.6}An epitrochoid is a plane curve traced  by a point
 of a circle which rolls without slipping on a fixed circle in the same plain, where the rolling circle rolls outside of the fixed circle.\end{Dn}
 Denote the radius of the fixed circle by $R$,  the radius of the rolling  circle by $r$ and the distance from the centre of the exterior circle   to the tracing
 point by $d$. Then, the parametric equations  for the epitrochoid have the form:
\begin{align}& x=\phi(t)=(R+r)\cos(t)-d\cos\left(\frac{R+r}{r}t\right),\label{8.11}
\\ &
 y=\psi(t)=(R+r)\sin(t)-d\sin\left(\frac{R+r}{r}t\right).\label{8.12}\end{align}
The epitrochoid \eqref{8.11}, \eqref{8.12} may be written in the form
\begin{equation}w=f(t)=(R+r)e^{it\lambda_1}-de^{it\lambda_2},\label{8.14}\end{equation}
where $\lambda_1=1,\,\, \lambda_2=(R+r)/r. $ If $(R+r)/r$ is an irrational number,
we deal with the case  discussed in 
Remark~\ref{Remark 8.3}. 
 
 If $(R+r)/r$ is a rational number, it may be represented in the form
 $(R+r)/r=p/q,$ where $p$ and $q$ are positive, coprime integers, $p>q.$ 
\begin{Pn}\label{Pn8.2+} System \eqref{8.11}, \eqref{8.12} is equivalent to the equation \eqref{5.1'}, where
\begin{equation} c_1=R+r,\quad c_2=-d,\quad m=q,\quad \ell=p, \quad t=qs.\label{8.13}\end{equation}\end{Pn}
The self-intersection points of the epitrochoid, where  $(R+r)/r$ is a rational number, are described
in Theorem \ref{Theorem 5.6} and in Corollary \ref{Corollary 5.8}.  
If $d=r$, then the epitrochoid coincides with the epicycloid.

We note that epitrochoids appear in the kinematics of the Wankel engine \cite{Lud, Nas}.

IV. \emph{Hypotrochoid}  \cite{BrG}.
\begin{Dn}\label{Definition 8.8}A hypotrochoid is a plane curve traced  by a point attached to a circle
which rolls without slipping on a fixed circle in the same plain, where the rolling circle rolls inside of the fixed circle.\end{Dn}
Denote the radius of the fixed circle by $R$,  the radius of the rolling  circle by $r$ and the distance from the centre of the interior circle   to the tracing
 point by $d$. 
 Then, the parametric equations
for the hypotrochoid have the form:
\begin{align}& x=\phi(t)=(R-r)\cos(t)+d\cos\left(\frac{R-r}{r}t\right),\label{8.15}
\\ &
y=\psi(t)=(R-r)\sin(t)-d\sin\left(\frac{R-r}{r}t\right).\label{8.16}\end{align}
The hypotrochoid \eqref{8.15}, \eqref{8.16}  may be written in the form
\begin{equation}w=f(t)=(R-r)e^{it\lambda_1}+de^{it\lambda_2},\label{8.18}\end{equation}
where $\lambda_1=1,\,\, \lambda_2=-(R-r)/r.$
If $(R-r)/r$ is an irrational number, we deal once more with the case  discussed in 
Remark~\ref{Remark 8.3}. 

 If $(R-r)/r$ is a rational number, it may be represented in the form
 $(R-r)/r=p/q,$ where $p$ and $q$ are positive,  coprime integers. 
Similar to the previous cases, we have the following proposition.
\begin{Pn}\label{Pn8.9} System \eqref{8.15}, \eqref{8.16} is equivalent to the equation \eqref{5.1'}, where
\begin{equation}  c_1=R-r,\quad  c_2=d\quad m=-q,\quad \ell=p, \quad  t=-qs.\label{8.17}\end{equation}\end{Pn}
The self-intersection points of the hypotrochoid, where  $(R+r)/r$ is a rational number, are described
in Theorem \ref{Theorem 5.6} and in Corollary \ref{Corollary 5.8}.  
If $d=r$, the hypotrochoid coincides with the hypocycloid.

\section{Curves in the three-dimensional space}
The parametric equations of a three-dimensional curve have the form
\begin{equation}x=\phi(t),\quad y=\psi(t),\quad z=\omega(t).\label{9.1}\end{equation}
We assume that $\phi(t), \,\psi(t),\, \omega(t)$ are smooth functions and $-\infty<t<+\infty$. The projection of a curve \eqref{9.1}
on the plane $(x,y)$ is a plane curve:
\begin{equation}x=\phi(t),\quad y=\psi(t).\label{9.2}\end{equation}
Let us describe the connections between the curves \eqref{9.1} and \eqref{9.2}. Self-intersection points of the curve \eqref{9.2}
appear when two different points of the curve \eqref{9.1} have the same projections on the plane $(x,y)$. Local singular points  of
the curve \eqref{9.2} appear when the tangent line to the curve \eqref{9.1} is parallel to the direction of projection, that is, to the axis $z$.
\begin{Rk}\label{Remark 9.1} Therefore, the methods of this paper may be useful not only in the theory of the plain curves, but in the
theory of the three-dimensional curves as well.\end{Rk}
We will need the following  definition (which is similar to Definition \ref{Definition 3.1} for the planar case).
\begin{Dn}\label{Definition 9.2} A curve \eqref{9.1} has  a self-intersection at some point $\break (x_0,y_0,z_0)$ if the curve 
has different tangent lines  at this point.\end{Dn}
Introduce the vectors
\begin{align}&
R(t):=(x(t),y(t),z(t))=(\phi(t),\psi(t),\quad \omega(t)).
\label{9.0}\end{align}
\begin{Pn}\label{Proposition 9.3} The curve \eqref{9.1} has different tangent lines at the point $(x_0,y_0,z_0)$  if and only if $R(t_1)=R(t_2)=(x_0,y_0,z_0)$
but the derivatives  $\dot{R}(t_1)$ and $\dot{R}(t_2)$ are not collinear  for some $t_1,t_2\in \BR$.\end{Pn}
In view of Definition \ref{Definition 3.1}, we have the corollary below.
\begin{Cy}\label{Corollary 9.4}If $R(t_1)=R(t_2)=(x_0,y_0,z_0)$ and if the point $(x_0,y_0)$ is a self-intersection point of the plain curve \eqref{9.2}, then
the point $(x_0,y_0,z_0)$ is a self-intersection point of the spatial curve \eqref{9.1}.\end{Cy}
\begin{Dn}\label{Definition 9.5}
The local singular points of the curve \eqref{9.1} are those points where \begin{equation}\frac{d\phi(t)}{dt}=\frac{d\psi(t)}{dt}=\frac{d\omega(t)}{dt}=0.\label{9.3}\end{equation}\end{Dn}
\begin{Cy}\label{Corollary 9.6}If $R(t_1)=(x_0,y_0,z_0),$ where  the point $(x_0,y_0)$ is a local singular point of the plain curve \eqref{9.2} and $\displaystyle{\frac{d\omega}{dt}(t_1)=0},$ then
the point $(x_0,y_0,z_0)$ is a local singular point  of the spatial curve \eqref{9.1}.\end{Cy}

\begin{Dn}\label{Definition 9.7} A $3$-dimensional curve $\Gamma$ is called a periodic helix   if the orthogonal projection of $\Gamma$
on the plain $(x,y)$ belongs to an $n$-member chain and $z=a\sin(Qt)$, where $a$ is a real number and $Q$ is a rational number.\end{Dn}
Some well-known $3$-dimensional curves belong to the periodic helix class of curves (e.g., Papareda class, satellite class, curves of the constant precession).
We consider a periodic helix  in Example \ref{Example 11.4} (see also Figure 4 in Appendix).
 \section{Triangular operators \\ with the difference kernels}
In this section, we continue our study of the operators with the difference kernels undertaken in  \cite{Sak1, Sak3, Sak4, Sak2}.
For this purpose, we use some of the obtained in the previous sections results.
Let us introduce the operator $J^{i\alpha}$ acting in $L^{2}(0,\omega)$:
\begin{equation}J^{i\alpha}f=\frac{1}{\Gamma(i\alpha+1)}\frac{d}{dx}\int_{0}^{x}f(t)(x-t)^{i\alpha}dt,
\label{10.1}\end{equation}
where $\alpha=\overline{\alpha}$ and $\Gamma(z)$ is Euler gamma  function. We proved the following assertion (see \cite{Sak1}).
\begin{Pn}\label{Propostion 10.1} The spectrum of the operator $J^{i\alpha}$ coincides with the set:
\begin{equation}e^{-\pi |\alpha|/2}{\leq}|z|{\leq}e^{\pi |\alpha|/2}.\label{10.2}\end{equation}\end{Pn}
Further we study the operators $S$ of the form
\begin{equation} S=\sum_{k=1}^{n}c_{k}J^{i\alpha_k},\label{10.3}\end{equation}
where the coefficients $c_k$ are real and $\alpha_k$ are rational numbers. We represent the numbers $\alpha_k$ in the form
$\alpha_k=\frac{m_k}{N}$, where $m_k$ are integers and $N$ is the smallest common denominator of $\alpha_k$ $(1\leq k\leq n)$. Hence, we have
\begin{equation} S=r(J^{i/N}),\label{10.4}\end{equation}
where $r(z)$ is a polynomial
\begin{equation} r(z)=\sum_{k=1}^{n}c_{k}z^{m_k}.\label{10.5}\end{equation}
The spectral set of the operator $S$ is denoted by $\sigma(S)$. 
According to  the well-known property of rational functions \cite{Lee}, we have
\begin{equation}\sigma(r(J^{i/N}))=r(\sigma(J^{i/N})).\label{10.6}\end{equation}
Recall that in view of Proposition \ref{Propostion 10.1}  (see \cite{Sak1, Sak3} for details) the boundary of $\sigma(J^{i/N})$ consists of two circumferences:
\begin{equation} z= e^{it}e^{\pi/(2N)} \quad {\mathrm{and}} \quad z= e^{it}e^{-\pi/(2N)}\quad (0{\leq}t{\leq}2\pi).\label{10.7}\end{equation}
Taking into account \eqref{10.3},\eqref{10.6} and \eqref{10.7} we obtain the assertion.
\begin{Tm}\label{Theorem 10.2}The boundary of $\sigma(S)$ consists of two parts:
\begin{equation}f_1(t)=\sum_{k=1}^{n}c_{k}e^{itm_k}e^{m_k\pi/(2N)},\quad f_2(t)=\sum_{k=1}^{n}c_{k}e^{itm_k}e^{-m_k\pi/(2N)}.\label{10.8}\end{equation}
\end{Tm}
\begin{Rk}\label{Remark 10.3}The functions $f_1(t)$ and $f_2(t)$ are $n$-member chains. Thus, the results of the present paper may be applied
to the  analysis of the structure  of $\sigma(S)$ $($see Example \ref{Example 11.5}$)$.\end{Rk}
\begin{Rk}\label{Remark 10.4}  The structure of $\sigma(S)$ in the case of   the numbers $\alpha_k$ linearly independent in the field of rational numbers
is studied in \cite{Sak2}.\end{Rk}
\section{Examples}
In this section,  we study several examples of the periodic functions          and calculate their self-intersection  points.
\begin{Ee}\label{Example 11.1} Consider a 2-member chain:
\begin{equation}w=f(t)=e^{2it}+ce^{3it},\quad c=\overline{c} \not=0.\label{11.1}\end{equation}\end{Ee}
The curve above has the following  parametric form:
\begin{align} & x=\cos(2t)+c\cos(3t),\quad y=\sin(2t)+c\sin(3t).\label{11.2}\end{align}

Since $\ell-m=1$, according to Theorem \ref{Theorem 5.4} and Corollary \ref{Corollary 5.8} all self-intersection points are real,
that is, $y(t)=0$ there.
In view of the first equality in \eqref{2.7} and the second equality in \eqref{11.2}, 
the points of the curve \eqref{11.1}, where $y=0$, satisfy the equation
\begin{equation} \sin(2t)+c\sin(3t)=\sin(t)[2\cos(t)-c+4c\cos^{2}(t)]=0.\label{11.3}\end{equation}

Clearly, \eqref{11.3} holds if $\sin(t)=0$. For the case $\sin(t)=0, \quad |c|\not=2/3$
we have $\dot{x}(t)=0$ and $\dot{y}(t)\not=0$. According to Remark \ref{Remark 3.6}, there
is no self-intersection in this case. The case $|c|=2/3$ is discussed separately in Conclusion 3) below.

Two further sets of the roots of   \eqref{11.3} are given by the relations
\begin{align}& \cos(t_1)=\frac{-1+\sqrt{1+4c^2}}{4c}=\frac{c}{1+\sqrt{1+4c^2}},\label{11.4}
\\ &
 \cos(t_2)=\frac{-1-\sqrt{1+4c^2}}{4c}.\label{11.5}\end{align}
Equality \eqref{11.4} implies that
\begin{equation} |\cos(t_1)|<\frac{1}{2},\label{11.6}\end{equation}
and so there exists $t_1$ satisfying \eqref{11.4}. 
On the other hand, the inequality $ \left|\frac{-1-\sqrt{1+4c^2}}{4c}\right|{\leq}1$
is valid if and only if
$|c|{\geq}2/3$
Hence, $t_2$ satisfying \eqref{11.5} exists if and only if  $|c|\geq 2/3$.
Moreover, if $t_1\,$ ($t_2$) satisfies \eqref{11.4} (satisfies \eqref{11.5}), then $s=-t_1$ ($s=-t_2$) satisfies
the same equality.
Finally, from the proposition below, which is proved  in Appendix \ref{Proof}, we see that relations \eqref{3.10n} hold
and we may apply Theorem \ref{Theorem 3.9}.
\begin{Pn}\label{Proposition 11.8} If relation \eqref{11.3} is valid, $\sin(t)\not=0$ and $|c|\not=2/3$, then
\begin{equation}\dot{x}(t)\not=0,\quad \dot{y}(t)\not=0.\label{11.18} \end{equation}\end{Pn}
We set $x_k=x(t_k)$ $(k=1,2)$, where $t_k$ are introduced in \eqref{11.4} and \eqref{11.5}. Now, two conclusions follow from Theorem \ref{Theorem 3.9}.\\
{\large{Conclusions 1) and 2).}}\\
1)  If $|c|<2/3$, the curve \eqref{11.1} has one and only one  self-intersection point, namely,  the self-intersection point  $(x_1,0)$.\\
2) If  $|c|>2/3$,  the curve \eqref{11.1} has  two and only two self-intersection points, namely,  the self-intersection points  $(x_1,0)$ and $(x_2,0)$.

If  $c=2/3$, taking into account \eqref{2.7}, \eqref{11.2}, \eqref{11.4} and \eqref{11.5}, we have
\begin{equation} \cos(t_1)=1/4,\quad x_1=-4/3,\quad \cos(t_2)=-1,\quad x_2=1/3,\label{11.9}\end{equation}
where $x_k=x(t_k)$.
If $c=-2/3$ , we obtain
\begin{equation} \cos(t_1)=-1/4,\quad x_1=-4/3,\quad \cos(t_2)=1,\quad x_2=1/3.\label{11.10}\end{equation}

The third equalities in \eqref{11.9} and \eqref{11.10} and formula \eqref{11.22} imply that $\sin(t_2)=0$ and $\dot{y}(t_2)=0$
for the case $|c|=2/3$.
Hence, from the proof of Proposition \ref{Proposition 11.8}  follows a corollary.
\begin{Cy}\label{CyProof} If $|c|=2/3$,  we have
\begin{equation} \dot{x}(t_1)\not=0, \quad \dot{y}(t_1)\not=0;\quad  \dot{x}(t_2)= \dot{y}(t_2)=0.\nonumber\end{equation}\end{Cy}
Now,  Corollary \ref{Corollary 5.8} and Theorems \ref{Theorem 7.3}   and \ref{Theorem 7.4} yield two more conclusions.\\
{\large{Conclusions 3) and 4).}}\\
3) If $|c|=2/3$,  the curve \eqref{11.1} has one self-intersection point and one local
singular point. These 
points are  $(x_1,0)$  and  $(x_2,0)$,  respectively, where $x_1,x_2$ for $c=2/3$ are given in \eqref{11.9} and $x_1,x_2$ for $c=-2/3$ are given in \eqref{11.10}.\\
4) The curve \eqref{11.1} has no   local singular points in the case $|c|\not=2/3.$\\

The plot of the curve \eqref{11.1}  (in the case $c=-2/3$)  illustrates the example (see Appendix \ref{Fig}, Figure 1).
\begin{Ee}\label{Example 11.2} Let us consider the parametric system
\begin{equation}\phi(t)=\cos(t)+\cos(6t),\quad \psi(t)=\sin(t)+\sin(6t),\label{11.11}\end{equation}\end{Ee}
Here, we have
\begin{equation}\ell=6,\quad m=1,\quad \ell-m=5.\label{11.12}\end{equation}
The plot of the curve \eqref{11.11}  (see Appendix \ref{Fig}, Figure 2) is an illustration of Corollary  \ref{Corollary 5.8}.
\begin{Ee}\label{Example 11.3} Let us consider the parametric system
\begin{equation}\phi(t)=\cos(t)+\cos(7t),\quad \psi(t)=\sin(t)+\sin(7t).\label{11.13}\end{equation}\end{Ee}
In the case \eqref{11.13}, we have
\begin{equation}\ell=7,\quad m=1,\quad \ell-m=6.\label{11.14}\end{equation}
The plot of the curve \eqref{11.13}  (see Appendix \ref{Fig}, Figure 3)  is an illustration of Corollary  \ref{Corollary 5.9}.
\begin{Ee}\label{Example 11.4} Let us consider a $3$-dimensional curve  $($see section 9$):$
\begin{equation}\phi(t)=\cos(t)+\cos(6t),\quad \psi(t)=\sin(t)+\sin(6t), \quad  \omega(t)=\sin(t). \label{11.15}\end{equation}
\end{Ee}
 Projection of the curve \eqref{11.15}
on the plane$(x,y)$ is a plane curve \eqref{11.11}. See Appendix \ref{Fig}, Figures 2 and 4, which   are the plots
of the plain curve \eqref{11.11}  and the space curve \eqref{11.15},  respectively.
We note that  Figure 4  is a periodic helix (see Definition 9.5).
\begin{Ee}\label{Example 11.5} Let us consider the operator
\begin{equation} S=2J^{i}+J^{2i},\label{11.16}\end{equation}
where the  operator $J^{i\alpha}$ is defined by relation \eqref{10.1}. \end{Ee}
In view of Theorem 10.2, we have the following result.
\begin{Tm}\label{Theorem 11.6}The boundary of $\sigma(S)$ consists of two parts:
\begin{equation}f_1(t)=\sum_{k=1}^{2}c_{k}e^{itm_k}e^{m_k\pi/2},\quad f_2(t)=\sum_{k=1}^{2}c_{k}e^{itm_k}e^{-m_k\pi/2},\label{11.17}\end{equation}
where $c_1=2,\quad c_2=1,\quad m_1=1,\quad m_2=2.$
\end{Tm}
\begin{Rk}\label{Remark 11.7.} The functions $f_1(t)$ and $f_2(t)$ are $2$-member chains.\end{Rk}
\section{Periodic helix}
Let us consider the periodic helix:
\begin{equation}x(t)=c_1\cos(mt)+c_2\cos({\ell}t),\nonumber\end{equation}
\begin{equation} y(t)=c_1\sin(mt)+c_2\sin({\ell}t),\,\label{12.1}\end{equation}
\begin{equation} z(t)=a\sin(\theta+Qt),\label{12.2}\end{equation}
where $c_1,c_2,m,\ell, a,\theta, Q$ are real and $\ell>|m|,\ \ell, and m$are integers.
In case of  periodic helix we have
\begin{equation}x^{2}(t)+y^{2}(t)=(c_1+c_2)^2-4c_1c_2\sin^{2}\frac{\ell-m}{2}t.\label{12.3}\end{equation}
The well-known
 Capareda curves, curves of constant precession \cite{Sco},\cite{LSTT}) belong to the class of periodic helix curves.
\begin{Pn}\label{Proposition 12.1}If
\begin{equation} c_1c_2>0,\theta=0,Q=\frac{\ell-m}{2},a^2=4c_1c_2\label{12.4}\end{equation} then
\begin{equation}x^2(t)+y^2(t)+z^2(t)=(c_1+c_2)^2.\label{12.5}\end{equation}\end{Pn}
In case \eqref{12.4} we have the  Capareda curves. Capareda curves are the curves traced on the sphere (see \eqref{12.5}).
\begin{Pn}\label{Proposition 12.2}If
\begin{equation}c_1c_2<0,\theta=0,Q=\frac{\ell-m}{2},a^2=-4c_1c_2\label{12.6}\end{equation} then
\begin{equation}x^2(t)+y^2(t)-z^2(t)=(c_1+c_2)^2.\label{12.7}\end{equation}\end{Pn}
In case \eqref{12.6} we have the   curves of constant precession.The  curves of constant precession
 are  traced on the hyperboloid (see \eqref{12.7}).
 It follows from Proposition 7.1 the assertion:
\begin{Pn}\label{Proposition 12.3}If
\begin{equation}|mc_1|{\ne}|\ell{c_2}| \label{12.8}\end{equation}then the corresponding  helix has  no singular points.\end{Pn}
\begin{Dn}\label{Definition 12.4}A periodic helix that  satisfies the conditions
\begin{equation}\theta=0,\quad Q=(\ell-m)/2,\label{12.9}\end{equation}
 we will call a S-periodic helix.\end{Dn}
\begin{Rk}\label{Remark 12.5} The Capareda curves and the   curves of constant precession belong to the class of
S-periodic helices.\end{Rk}
Using relations \eqref{12.2}, \eqref{12.3} and \eqref{12.9}we obtain:
\begin{Tm}\label{Theorem 12.6}In case of S-periodic helix the following assertions are valid:\\
1. If $c_1c_2>0$  then
\begin{equation}x^2(t)+y^2(t)+z^2(t)=(c_1+c_2)^2+(a^2- 4c_1c_2)\sin^{2}\frac{\ell-m}{2}t.\label{12.10}\end{equation}
2. If $c_1c_2<0$  then
\begin{equation}x^2(t)+y^2(t)-z^2(t)=(c_1+c_2)^2-(a^2+ 4c_1c_2)\sin^{2}\frac{\ell-m}{2}t.\label{12.11}\end{equation}\end{Tm}
 The Theorem 12.6 allow us to make conclusions about the position of  the corresponding S-periodic helix relative to the sphere
(assertion 1) or relative to the hyperboloid (assertion 2).
\begin{Pn}\label{Proposition 12.7}In case of S-periodic helix the following assertions are valid:\\
1. If
\begin{equation} c_1c_2>0,\quad (a^2- 4c_1c_2){\geq}0, \label{12.12}\end{equation}
then
\begin{equation}(c_1+c_2)^2{\leq}x^2(t)+y^2(t)+z^2(t){\leq}(c_1+c_2)^2+(a^2- 4c_1c_2),\label{12.13}\end{equation}

2. If
\begin{equation} c_1c_2>0,\quad (a^2- 4c_1c_2){\leq}0, \label{12.14}\end{equation}
then
\begin{equation}(c_1+c_2)^2+(a^2- 4c_1c_2){\leq}x^2(t)+y^2(t)+z^2(t){\leq}(c_1+c_2)^2,\label{12.15}\end{equation}
3. If
\begin{equation} c_1c_2<0,\quad (a^2+ 4c_1c_2){\leq}0, \label{12.16}\end{equation}
then
\begin{equation}(c_1+c_2)^2+(a^2+ 4c_1c_2){\leq}x^2(t)+y^2(t)-z^2(t){\leq}(c_1+c_2)^2,\label{12.17}\end{equation}
4. If
\begin{equation} c_1c_2>0,\quad (a^2+ 4c_1c_2){\geq}0, \label{12.18}\end{equation}
then
\begin{equation}(c_1+c_2)^2{\leq}x^2(t)+y^2(t)-z^2(t){\leq}(c_1+c_2)^2+(a^2+ 4c_1c_2),\label{12.19}\end{equation}
\end{Pn}

\section{Torus Knots}
\begin{Dn}Knot is a three-dimensional curve  which is homeomorphic to a circle.\end{Dn}
The equations of torus knots can be written in the parametric form:\cite{Kaw}
\begin{equation}x(t)=[R+r\cos(qt)]\cos(pt),\label{13.1}\end{equation}
\begin{equation}y(t)=[R+r\cos(qt)]\sin(pt),\label{13.2}\end{equation}
\begin{equation}z(t)=r\sin(qt),\label{13.3}\end{equation}
where $p$ and $q$ are positive integers and
\begin{equation}R>r>0.\label{13.4}\end{equation}
Every torus  knot can be expressed as Fourier -(3,3,1) knot.
Indeed,we may rewrite \eqref{13.1}-\eqref{13.3} in the following form
\begin{equation}x(t)=R\cos(pt)+(r/2)\cos(p+q)t+(r/2)\cos(p-q)t,\label{13.5}\end{equation}
\begin{equation}y(t)=R\sin(pt)+(r/2)\sin(p+q)t-(r/2)\sin(q-p)t,\label{13.6}\end{equation}
\begin{equation}z(t)=r\sin(qt).\label{13.7}\end{equation}
It follows from \eqref{13.5} and \eqref{13.6} that
\begin{equation}[R-\sqrt{x^{2}(t)+y^{2}(t)}]^{2}=r^2-r^2\sin^{2}(qt)\label{13.8}\end{equation}
For torus knots we have
\begin{equation}[R-\sqrt{x^{2}(t)+y^{2}(t)}]^{2}+z^{2}(t)=r^2\label{13.9}\end{equation}
Relation \eqref{13.9}is equation of torus.Hence the knot torus curves are the curves traced on the torus,
\begin{Dn}\label{Definition 13.2}A periodic helix that satisfies relations \eqref{13.5}. \eqref{3.6} and the equality
\begin{equation}z(t)=a\sin(qt),\quad a>0\label{13.10}\end{equation}
we will call S-torus knot.\end{Dn}
\begin{Rk}\label{Remark 13.3}The torus knots belong to the class S-torus knots.\end{Rk}
\begin{Tm}\label{Theorem 13.4}In case of S-torus knot the following relation is valid
\begin{equation}[R-\sqrt{x^{2}(t)+y^{2}(t)}]^{2}+z^{2}(t)=r^2+(a^2-r^2)\sin^{2}(qt)\label{13.11}\end{equation}\end{Tm}
The Theorem 13.4 allow us to make conclusions about the position of  the corresponding S-periodic helix relative to the torus.
\begin{equation}[R-\sqrt{x^{2}(t)+y^{2}(t)}]^{2}+z^{2}(t)=r^2+(a^2-r^2)\sin^{2}(qt)\label{13.12}\end{equation}
\begin{Pn}\label{Proposition 13.5}In case of S-torus knot the following assertions are valid:\\
1. If $a{\geq}r$ then
\begin{equation}0{\leq}[R-\sqrt{x^{2}(t)+y^{2}(t)}]^{2}+z^{2}(t)-r^2{\leq}(a^2-r^2)\label{13.13}\end{equation}
2. If $a{\leq}r$ then
\begin{equation}(a^2-r^2){\leq}[R-\sqrt{x^{2}(t)+y^{2}(t)}]^{2}+z^{2}(t)-r^2{\leq}0\label{13.14}\end{equation}
\end{Pn}
\section{Self-intersection points}In this section we investigate the self-intersection points of the plane curves \eqref{13.1},\eqref{13.2}( the orthogonal
projections of the torus knots onto the plane (x,y)).\\
Using  the trigonometrical formulas
\begin{equation}
\sin(pt)-\sin(p\tau)=2\cos[p(t+\tau)/2]\sin[p(t-\tau)/2].\label{14.1}\end{equation}
\begin{equation}
\cos(pt)-\cos(p\tau)=-2\sin[p(t+\tau)/2]\sin[p(t-\tau)/2].\label{14.2}\end{equation}
we obtain the assertion:
\begin{Pn}\label{Proposition 14.1}
If $t{\ne}\tau , (mod\,(2\pi)$ and
\begin{equation}\sin(pt)=\sin(p\tau),\,\cos(pt)=\cos(p\tau),\,\cos(qt)=\cos(q\tau) \label{14.3}\end{equation}
then
\begin{equation}t_{k,j}-\tau_{k,j}=2k\pi/p \quad (mod \,2\pi),\label{14.4}\end{equation}
\begin{equation}t_{k,j}+\tau_{k,j}=2j\pi/q \quad (mod \,2\pi).\label{14.5}\end{equation}
where $ p{\geq}2,\,1{\leq}k{\leq}p-1,\quad  0{\leq}j{\leq}q-1.$\end{Pn}
Hence we have
\begin{equation}t_{k,j}=\frac{k\pi}{p}+\frac{j\pi}{q},\quad  \tau_{k,j}=-\frac{k\pi}{p}+\frac{j\pi}{q},\quad (mod \, 2\pi),\label{14.6}\end{equation}
where $1{\leq}k{\leq}p-1,\quad  0{\leq}j{\leq}q-1.$
\begin{Tm}\label{Theorem 14.2} Let $x(t)$ and $y(t)$ be defined by relations \eqref{13.1} and \eqref{13.2}.The relations
\begin{equation}
x(t)=x(\tau),\quad  y(t)=y(\tau),\quad  (t{\ne}\tau,\mod\,2\pi)\label{14.7}\end{equation}
are valid if and only if  the equalities \eqref{4.6} are valid.\end{Tm}
\emph{Proof.} It follows from (12.1) and (14.7) that
\begin{equation}[R+r\cos(qt)]^2=[R+r\cos(q\tau)]^2.\label{14.8}\end{equation}
In the case
\begin{equation}[R+r\cos(qt)]=[R+r\cos(q\tau)].\label{14.9}\end{equation}
the conditions \eqref{14.3} are fulfilled. Then the Theorem  follows from Proposition 14.1.
 Let us consider the case when
\begin{equation}[R+r\cos(qt)]=-[R+r\cos(q\tau)].\label{14.10}\end{equation}
Hence we have
\begin{equation}2R+r[\cos(qt)+\cos(q\tau)]=0.\label{14.11}\end{equation} The  relation \eqref{14.11} is impossible because $R>r$. Then the case (14.12) is impossible too.
So we have proved the Theorem.\\
Let us now prove that the points \eqref{14.6} are  self-intersection points of the curve \eqref{13.1}, \eqref{13.2} (see Definition 3.1 and Proposition 3.2).
\begin{La}\label{Lemma 14.3}The following inequalities
\begin{equation}\sin(qt_{k,j}){\ne}\sin(q\tau_{k.j}),\quad 1{\leq}k{\leq}p-1,\quad 0{\leq}j{\leq}q-1.\label{14.12}\end{equation}
are valid.\end{La}
\emph{Proof.} Using relations  (14.4)-(14.5)    and formula
 \begin{equation}
\sin(qt_{k,j})-\sin(q\tau_{k,j})=\nonumber\end{equation}
\begin{equation}2\cos[q(t_{k,j}+\tau_{k,j})/2]\sin[q(t_{k,j}-\tau_{k,j})/2]],\label{14.13}\end{equation}
we have
\begin{equation}
\sin(qt_{k,j})-\sin(q\tau_{k.j})=\nonumber\end{equation}
\begin{equation}2\sin(qk\pi/p), \quad 1{\leq}k{\leq}p-1,\quad 0{\leq}j{\leq}q-1.\label{14.14}\end{equation}
Since $k$ is less than $p$, and the numbers $p$ and $q$ are coprime, then the number $qk/p$ is not an integer.Hence
\begin{equation}\sin(qk\pi/p){\ne}0, \quad 1{\leq}k{\leq}p-1.\label{14.15}\end{equation} The Lemma is proved.\\
 Let us write
\begin{equation}\dot{x}(t)=-rq\sin(qt)\cos(pt)-[R+r\cos(qt)]p\sin(pt),\label{14.16}\end{equation}
\begin{equation}\dot{y}(t)=-rq\sin(qt)\sin(pt)+[R+r\cos(qt)]p\cos(pt),\label{14.17}\end{equation}
Relations (14.6), (14.12), (14.16) and (14.17) imply:
\begin{La}\label{Lemma 14.4} The following relations are valid
\begin{equation} \dot{x}(t_{k,0})=-rq\sin(qt_{k,0})\cos(pt_{k,0}){\ne}0,\quad  1{\leq}k{\leq}p-1,\label{14.18}\end{equation}
\begin{equation} \dot{x}(\tau_{k,0})=-rq\sin(q\tau_{k,0})\cos(p\tau_{k,0}){\ne}0,\quad  1{\leq}k{\leq}p-1,\label{14.19}\end{equation}
\begin{equation}\dot{y}(t_{k,0})=[R+r\cos(qt_{k,0})]p\cos(pt_{k,0})=\nonumber\end{equation}
\begin{equation}\dot{y}(\tau_{k,0}){\ne}0,\quad  1{\leq}k{\leq}p-1,\label{14.20}.\end{equation}
\begin{equation} \dot{y}(t_{k,0})/ \dot{x}(t_{k,0}){\ne} \dot{y}(\tau_{k,0})/ \dot{x}(\tau_{k,0}),\quad  1{\leq}k{\leq}p-1.
\label{14.21}
\end{equation}
\end{La}

Let us introduce the function (see (13.1) and (13.2)):
\begin{equation}f(t)=x(t)+iy(t)=e^{ipt}[R+r\cos(qt)].\label{14.22}\end{equation}
Then we have
\begin{equation}f(t+2j\pi/q))=e^{i\beta}f(t),\quad \dot{f}(t+2j\pi/q))=e^{i\beta}\dot{f}(t),\label{4.23}\end{equation}
\begin{equation}\dot{x}(t_0+2j\pi/q)=\cos(\beta)\dot{x}(t_0)-\sin(\beta)\dot{y}(t_0),\label{4.24}\end{equation}
\begin{equation}\dot{y}(t_0+2j\pi/q)=\sin(\beta)\dot{x}(t_0)+\cos(\beta)\dot{y}(t_0),\label{4.25}\end{equation}
where $\beta=2jp\pi/q.$ Using (14.6), (14.26) and (14.27) we obtain the following relation
\begin{equation}\dot{y}(t_{k,j})\dot{x}(\tau_{k,j})-\dot{y}(\tau_{k,j})\dot{x}(t_{k,j})= \nonumber\end{equation}
\begin{equation}\dot{y}(t_{k,0})\dot{x}(\tau_{k,0})-\dot{y}(\tau_{k,0})\dot{x}(t_{k,0}).\label{14.26}\end{equation}
From (14.21) and (14.26) we get
\begin{equation}\dot{y}(t_{k,j})\dot{x}(\tau_{k,j})-\dot{y}(\tau_{k,j})\dot{x}(t_{k,j}){\ne}0,1{\leq}k{\leq}p-1.\label{14.27}\end{equation}
Hence.
\begin{equation}\dot{y}(t_{k,j})/\dot{x}(t_{k,j}){\ne}\dot{y}(\tau_{k,j})/\dot{x}(\tau_{k,j})\label{14.28}\end{equation}
So we proved the following assertion
\begin{Tm}\label{Theorem 14.5}The points
\begin{equation}M_{k,j}[x(t_{k,j}),y(t_{k,j})],(1{\leq}k{\leq}p-1,\quad 0{\leq}j{\leq}q-1)\label{14.29}\end{equation}
are self-intersection points of the curve (13.1), (13.2).\end{Tm}
\begin{Cy}\label{Corollary 14.6}The curve (13.1),(13.2) have $qp-q$ self-intersection points.\end{Cy}
\emph{Proof.} Let us suppose that\begin{equation}\frac{k_1}{p}+\frac{j_1}{q}=\frac{k_2}{p}+\frac{j_2}{q}.\label{14.30}\end{equation}
Then we have
\begin{equation}(k_1-k_2)q=(j_2-j_1)p.\label{14.31}\end{equation}
Numbers $q$ and $p$ are coprime, $|j_2-j_1|<q$. Hence
\begin{equation}k_1-k_2=j_2-j_1=0.\label{14.32}\end{equation}
There are no two identical points in the set of self-intersection points $M_{k,j},$ where $(1{\leq}k{\leq}p-1,\quad 0{\leq}j{\leq}q-1)$
The corollary is proved.
\begin{Rk}\label{Remark 14.7}The self-intersection points of curve (13.1), (13,2) were studied in an interesting paper \cite{HoPi} The authors, choosing a specific curve (torus knot (p=3,q=7)) as a model
conducted a heuristic analysis of the situation. In our work, we offer a rigorous mathematical analysis of the situation,
based on the classical definition of self-intersection points.\end{Rk}
\begin{Pn}\label{Pn14.8} S-torus knots  do not have singular points.\end{Pn}
\emph{Proof.} If S-torus knot have a singular point then (see (13.1),(13.2)):
\begin{equation}\dot{x}(t)=-qr\sin(qt)\cos(pt)-\nonumber\end{equation}
\begin{equation}(R+r\cos(qt))p\sin(pt)=0,\label{14.33}\end{equation}
\begin{equation}\dot{y}(t)=-qr\sin(qt)\sin(pt)+\nonumber\end{equation}
\begin{equation}(R+r\cos(qt))p\cos(pt)=0,\label{14.34}\end{equation}
It follows from (14.33) and (14.34)that
\begin{equation}\sin(qt)=0,\quad (R+r\cos(qt))=0.\label{14.35}\end{equation}
The relation $(R+r\cos(qt))=0$ is impossible because $R>r$.
So we have proved the Proposition.
\section{Fold points, N=2}
1. Let us consider the curve \eqref{1.1} when $n=2$.Thus
\begin{equation} f(t)={c_1}e^{itm}+{c_2}e^{it\ell},\quad  \ell>m,\quad \ell{\ne}-m,\label{15.1}\end{equation}
where $m$ and $\ell$ are coprime integers.\\
\textbf{Definition}\\
The curve \eqref{15.1} has in the point $t$ the fold if either
\begin{equation}\dot{\phi}(t)=-m\sin(tm)-{\ell}c_{2}\sin(t\ell)=0,\quad \dot{\psi}(t){\ne}0\label{15.2}\end{equation}
or
\begin{equation}\dot{\psi}(t)=mc_{1}\cos(tm)+{\ell}c_{2}\cos(t\ell)=0,\quad \dot{\phi}(t){\ne}0. \label{15.3}\end{equation}
We suppose that
\begin{equation} (mc_1)^2=(\ell{c_2})^2.\label{15.4}\end{equation}

\begin{Pn}\label{Proposition 15.2} All local fold points of the  curve \eqref{15.1} are double, i.e.
\begin{equation}[\ddot{\phi}(t)]^2+[\ddot{\psi}(t)]^2{\ne}0.\label{15.5}\end{equation}\end{Pn}
(see proof of Proposition 7.2).\\
2. Now we consider separately the case when
\begin{equation}m{c_1}=\ell{c_2}.\label{15.6}\end{equation}
Taking into account  \eqref{15.6}  we have
\begin{equation}\dot{\phi}(t)=-2m{c_1}\sin(\frac{\ell+m}{2}t)\cos(\frac{\ell-m}{2}t) \nonumber \end{equation}
\begin{equation}\dot{\psi}(t)=2m{c_1}\cos(\frac{\ell+m}{2}t)\cos(\frac{\ell-m}{2}t) \label{15.7} \end{equation}
It follows from \eqref{15.2} and \eqref{15.7} that
\begin{equation}\sin(\frac{\ell+m}{2}t)=0,\label{15.8}\end{equation}
It follows from  \eqref{15.3} and \eqref{15.7}  that
\begin{equation} \cos(\frac{\ell-m}{2}t)=0. \label{15.9} \end{equation}if relation \eqref{15.3} is fulfilled.
Hence we proved the assertion.
\begin{Tm}\label{Theorem 15.3}Let condition \eqref{15.6} be fulfilled. Then\\
1)The relations \eqref{15.2} are valid if and only if
\begin{equation}t_k=\frac{2k}{\ell+m}\pi,\quad 0{\leq}k{\leq}\ell+m-1.\label{15.10}\end{equation}
2)The relations \eqref{15.3} are valid if and only if
\begin{equation}\tau_k=\frac{2k+1}{\ell+m}\pi,\quad 0{\leq}k{\leq}\ell+m-1.\label{15.11}\end{equation}
3) The set of fold points of the curve \eqref{15.1} consists from two sets  $w_k=f(t_k)$ and   $v_k=f(\tau_k).$\end{Tm}
3. Now we consider the case
\begin{equation}m{c_1}=-\ell{c_2}.\label{15.12}\end{equation}
Taking into account \eqref{15.12}  we have
\begin{equation}\dot{\phi}(t)=2m{c_1}\sin(\frac{\ell-m}{2}t)\cos(\frac{\ell+m}{2}t) \nonumber \end{equation}
\begin{equation}\dot{\psi}(t)=2m{c_1}\sin(\frac{\ell+m}{2}t)\sin(\frac{\ell-m}{2}t) \label{15.13} \end{equation}
Relations \eqref{15.13} are equivalent to the equality
\begin{equation}
\cos(\frac{\ell+m}{2}t)=0. \label{15.14} \end{equation} if relation \eqref{15.2} is fulfilled,
and are equivalent to the equality
\begin{equation} \sin(\frac{\ell+m}{2}t)=0. \label{15.15} \end{equation}if relation \eqref{15.3} is fulfilled.
Hence we proved the assertion.
\begin{Tm}\label{Theorem 15.4}Let condition \eqref{15.12} be fulfilled. Then\\
1)The relations \eqref{15.2} are valid if and only if
\begin{equation}\tau_k=\frac{2k}{\ell+m}\pi,\quad 0{\leq}k{\leq}\ell+m-1.\label{15.16}\end{equation}
2)The relations \eqref{15.3} are valid if and only if
\begin{equation}t_k=\frac{2k+1}{\ell+m}\pi,\quad 0{\leq}k{\leq}\ell+m-1.\label{15.17}\end{equation}
3) The set of fold points of the curve \eqref{15.1} consists from two sets  $w_k=f(t_k)$ and   $v_k=f(\tau_k).$
\end{Tm}
It follows from Theorems 15.3 and 15.4  the  assertion
\begin{Cy}\label{Corollary 15.5} Let the condition \eqref{15.4} be fulfilled. Then
\begin{equation}w_k=w_0e^{i2k\pi/(\ell+m)},\, v_k=v_0e^{i2k\pi/(\ell+m)},\, 0{\leq}k{\leq}\ell+m-1.\label{15.18}\end{equation}
\end{Cy}
 Let us consider the finite group of  rotations of complex plane  with multiplication operators
\begin{equation}g_k=e^{i2k\pi/(\ell-m)},\quad 0{\leq}k{\leq}\ell+m-1.\label{15.19}\end{equation}
We denote this group by $G(\ell,m).$ From relations \eqref{15.18} and \eqref{15.19} we have
\begin{Tm}\label{Theorem 15.6} Let the condition \eqref{15.4} be fulfilled.The elements of group $G(\ell,m)$ map the corresponding
 class of fold points  onto itself.\end{Tm}
\begin{Rk}\label{Remark 15.7}It is well-known that the groups of rotations of complex plane  are Lie groups. Then group $G(\ell,m)$ is Lie group.\end{Rk}
\begin{Rk}\label{Remark 15.8}The self-intersection and singular points of the periodic functions were investigated in section 1-11 of the paper.\end{Rk}
\appendix
\section{Proof of Proposition \ref{Proposition 11.8}}\label{Proof}
 We prove the proposition by contradiction. Suppose that
\begin{equation}\dot{x}(t)=\sin(t)[12c\cos^{2}(t)+4\cos(t)-3c]=0.\label{11.19} \end{equation}
It follows from  \eqref{11.3} and \eqref{11.19} that $\cos(t)=0$. Using this fact and \eqref{11.3} we obtain the equality $c=0$.
From the last relation and \eqref{11.1} we have $\dot{x}(t)\not=0$ (i.e., we arrive at a contradiction).
The first inequality in Proposition \ref{Proposition 11.8} is proved.

Next,  note that
\begin{align}\label{11.20}  \dot{y}(t)&=12c\cos^{3}(t)+4\cos^{2}(t)-9c\cos(t)-2
\\ &
=3\cos(t)\big(4c\cos^2(t)+2\cos(t)-c\big)-2\big(\cos^2(t)+3c\cos(t)+1\big)
\nn \end{align}
Since $\sin(t)\not= 0$, relation \eqref{11.3} implies that 
\begin{align}\label{11.20!}  
& 4c\cos^2(t)+2\cos(t)-c=0,
\end{align}
and we rewrite \eqref{11.20} as
\begin{align}\label{11.20+} & \dot{y}(t)=-2g(t,c),  \quad g(t,c):=\cos^{2}(t)+3c\cos(t)+1.
\end{align}
The solutions of \eqref{11.20!} consist of the solutions $t_1$ given by \eqref{11.4} and solutions $t_2$ given by \eqref{11.5}.
According to \eqref{11.4}, the solutions $t_1$    satisfy the inequality
$c\cos(t)>0$.  It follows that  $g(t,c)>0$ or, equivalently,  $\dot{y}(t)<0$. 
The second inequality in Proposition \ref{Proposition 11.8} is proved for $t_1$.

The solutions $t_2$ of \eqref{11.20!} exist only in the case $|c|\geq 2/3$.
Taking into account that $|c|\not=2/3$
in Proposition~\ref{Proposition 11.8}, we study the case $|c|> 2/3$.
Suppose that
\begin{align}\label{11.21} & g(t,c)=0, \quad{\mathrm{that\,\, is,}} \quad \dot{y}(t)=0.\end{align}
Using \eqref{11.20!}--\eqref{11.21} we have
\begin{equation} \cos(t)=-\frac{3c^2+2}{5c}.\label{11.23}\end{equation}
Hence, the following equality is valid:
\begin{equation}-\frac{3c^2+2}{5c}=\frac{-1-\sqrt{1+4c^2}}{4c}.\label{11.24}\end{equation}
It follows from \eqref{11.24} that
\begin{equation}\sqrt{1+4c^2}=\frac{12c^2+3}{20}.\label{11.25}\end{equation}
Relation \eqref{11.24} implies that
\begin{equation}144c^4-1528c^2-391=0.\label{11.26}\end{equation}
Equation \eqref{11.26} is a quadratic equation with respect to $c^2$, which has two solutions $c_1^2=-1/4$ and $c_2^2=391/36$.
Clearly, the equality $c_1^2=-1/4$ does not hold for the real-valued $c_1$. For the case $c^2=391/36$, relation \eqref{11.23}
implies that $|\cos(t)|>1$ and we arrive at a contradiction (with the assumption $g(t,c)=0$).
Thus, the second inequality in Proposition \ref{Proposition 11.8} is proved for $t_2$ as well.
In this way, the Proposition \ref{Proposition 11.8} is proved.

Let us rewrite \eqref{11.20+} for the cases $c=2/3$ and $c=-2/3$:
\begin{equation}\dot{y}(t,2/3)=-2[\cos(t)+1]^2=0,\quad \dot{y}(t,-2/3)=-2[\cos(t)-1]^2=0.\label{11.22}\end{equation}

\newpage
\section{Figures}\label{Fig}
\begin{figure}[h]
    \centering
    \includegraphics[width=0.42\textwidth]{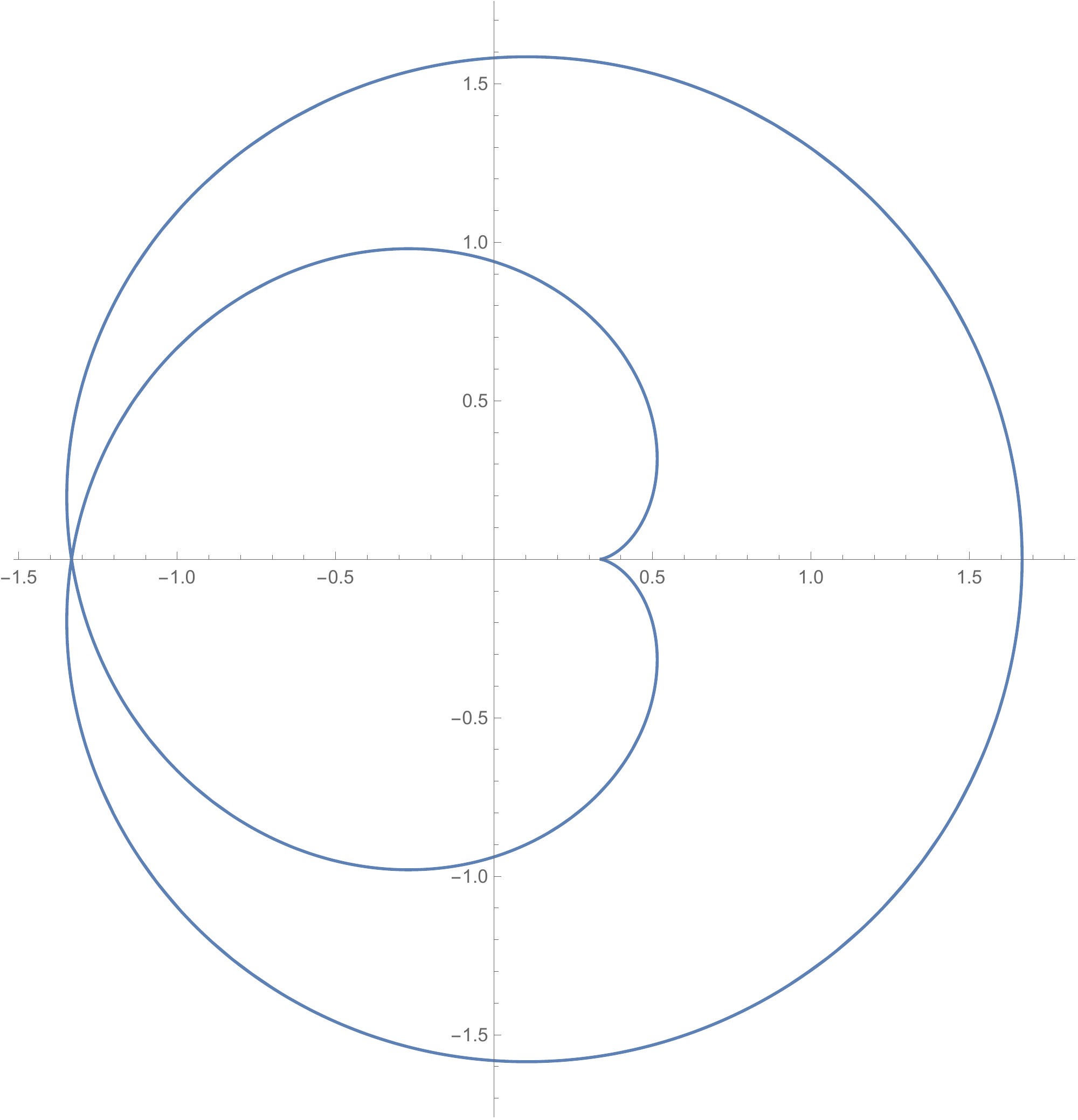}
    \caption{Example 11.1 ($c = -\frac{2}{3}$)}
    \label{fig:fig1}
\end{figure}

\begin{figure}[h]
    \centering
    \includegraphics[width=0.42\textwidth]{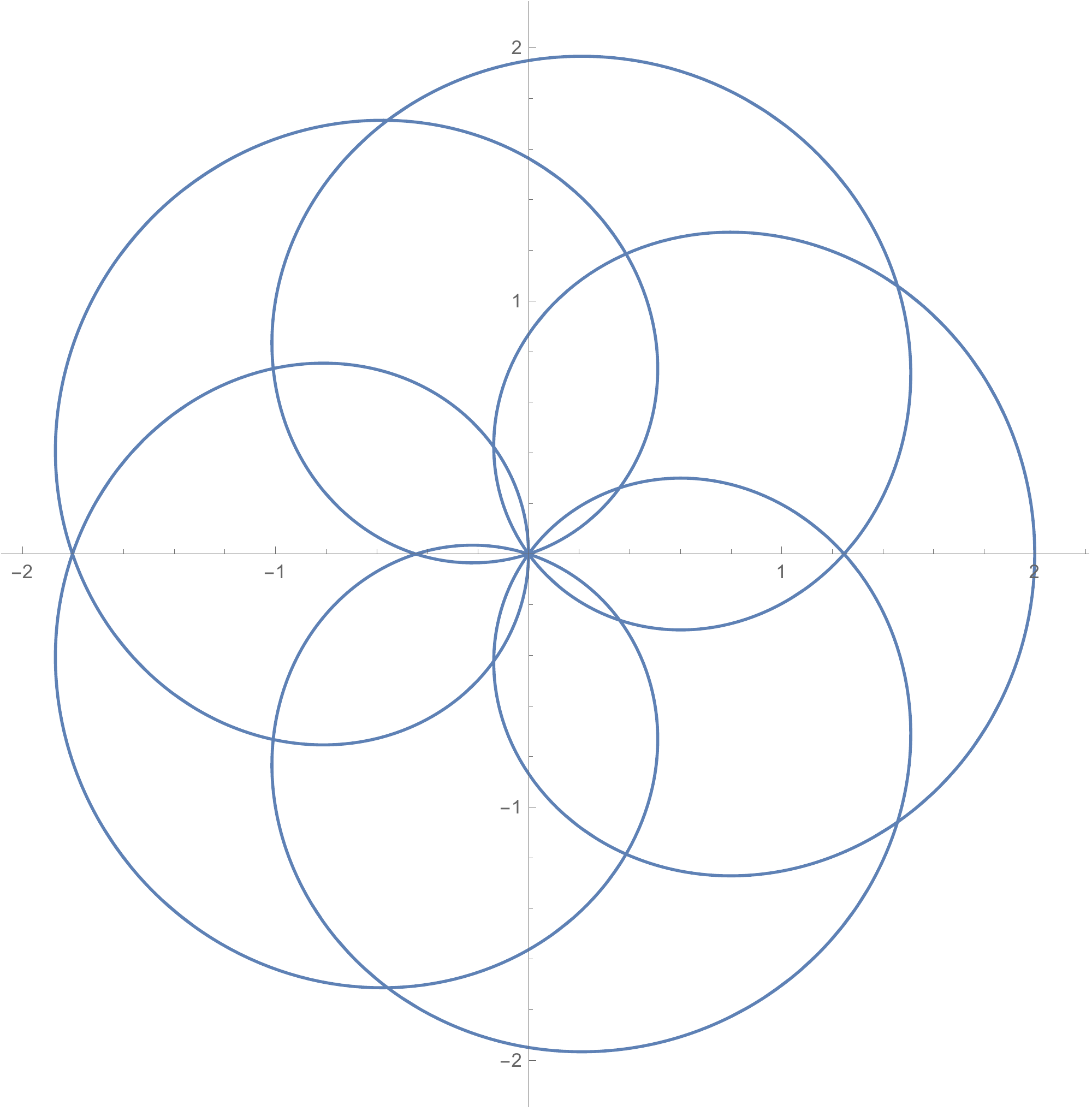}
    \caption{Example 11.4}
    \label{fig:fig2}
\end{figure}

\begin{figure}[h]
    \centering
    \includegraphics[width=0.42\textwidth]{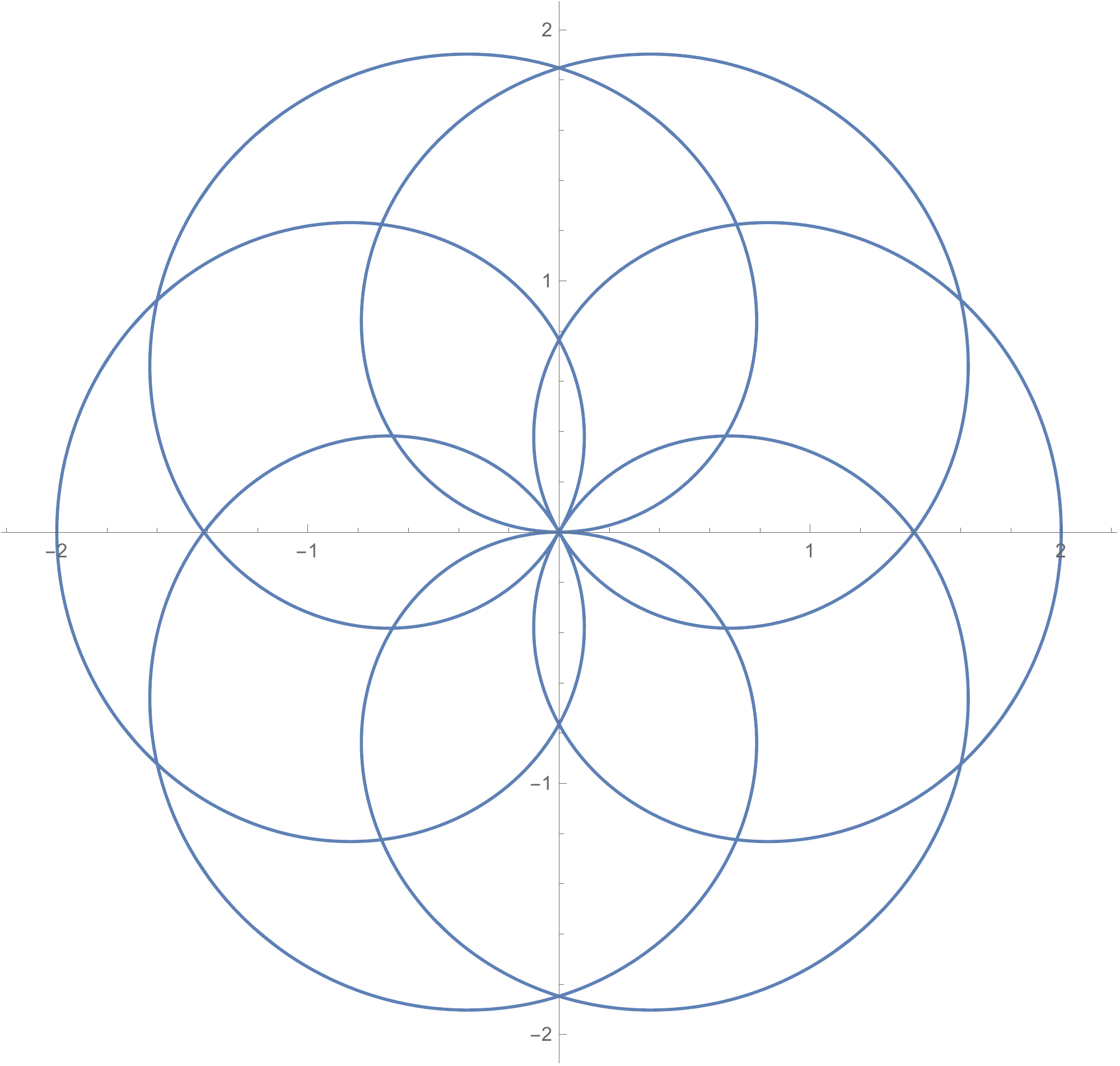}
    \caption{Example 11.5}
    \label{fig:fig3}
\end{figure}

\begin{figure}[h]
    \centering
    \includegraphics[width=0.45\textwidth]{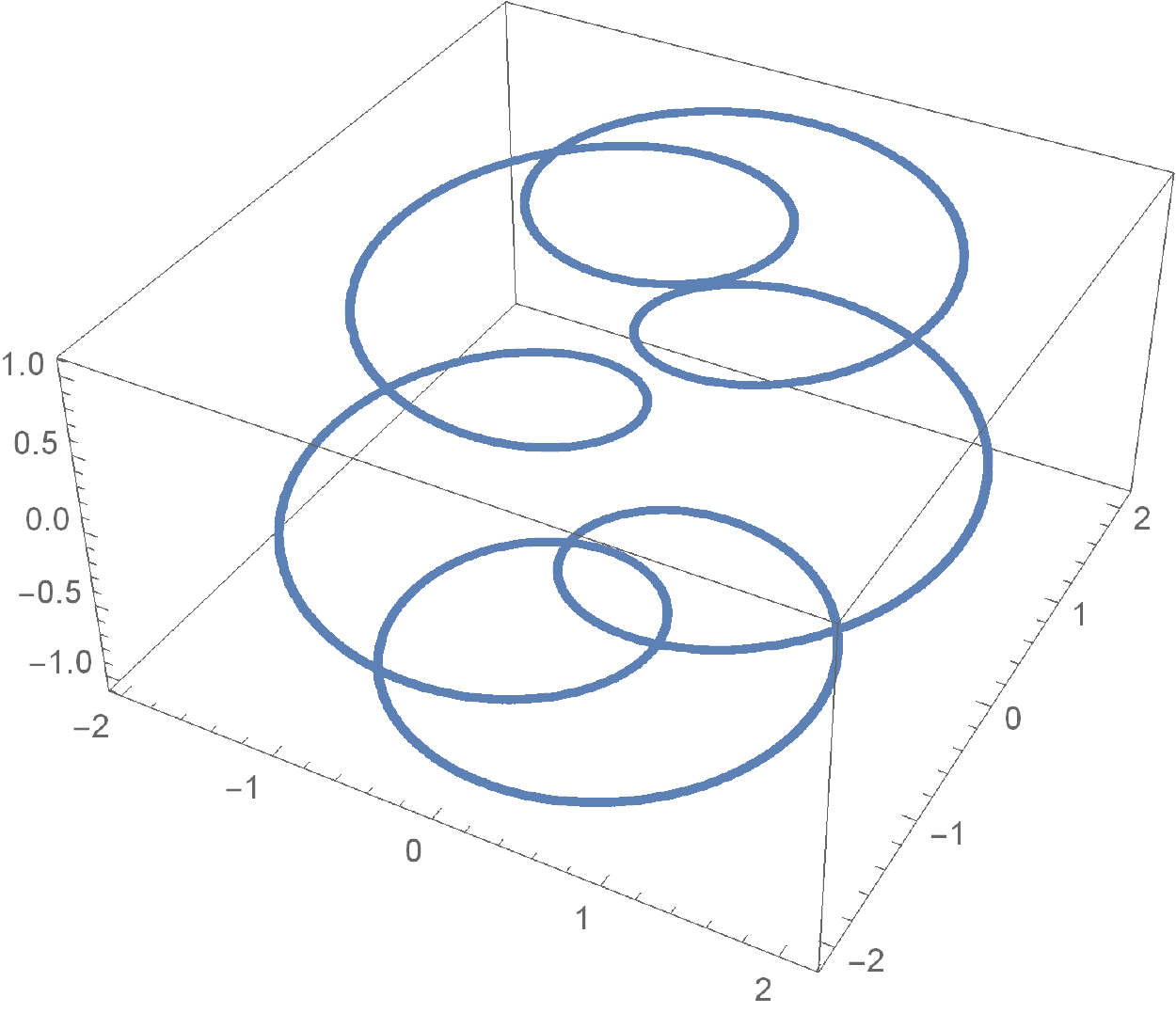}
    \caption{Example 11.6}
    \label{fig:fig4}
\end{figure}

\newpage

{\bf Acknowledgments}  {The author is grateful to A.L. Sakhnovich for useful discussions and important remarks.} 

\end{document}